\documentclass[11pt, letterpaper]{amsart}
\usepackage[margin=1.5in]{geometry}

\usepackage[english]{babel}
\usepackage[utf8]{inputenc}
\usepackage[T1]{fontenc}
\usepackage{lmodern} 

\usepackage[colorinlistoftodos]{todonotes}

\usepackage{amsfonts} 
\usepackage{amssymb}
\usepackage{amsmath}
\usepackage{amsthm}
\usepackage{graphicx}
\usepackage{xifthen}
\usepackage{verbatim}
\usepackage{enumerate}
\usepackage{needspace}
\usepackage{mathtools}
\usepackage{xfrac}
\usepackage{tikz-cd}
\usepackage{lipsum}

\usepackage[colorlinks=true, allcolors=blue]{hyperref}

\usepackage[nameinlink]{cleveref} 

\DeclareRobustCommand{\SkipTocEntry}[5]{}

\newcommand\blfootnote[1]{%
	\begingroup
	\renewcommand\thefootnote{}\footnote{#1}%
	\addtocounter{footnote}{-1}%
	\endgroup
}

\newtheorem*{lemma*}{Lemma}

\newtheorem{theorem}{Theorem}[section]
\newtheorem{proposition}[theorem]{Proposition}
\newtheorem{corollary}[theorem]{Corollary}
\newtheorem{lemma}[theorem]{Lemma}

\theoremstyle{definition}
\newtheorem{definition}[theorem]{Definition}

\theoremstyle{definition}
\newtheorem{example}[theorem]{Example}

\theoremstyle{definition}
\newtheorem{remark}[theorem]{Remark}

\crefname{chapter}{chapter}{chapters}
\Crefname{chapter}{Chapter}{Chapters}
\crefname{section}{section}{sections}
\Crefname{section}{Section}{Sections}
\crefname{subsection}{section}{sections}
\Crefname{subsection}{Section}{Sections}
\crefname{subsubsection}{section}{sections}
\Crefname{subsubsection}{Section}{Sections}
\crefname{figure}{figure}{figures}
\Crefname{figure}{Figure}{Figures}
\crefname{table}{table}{tables}
\Crefname{table}{Table}{Tables}

\crefname{theorem}{theorem}{theorems}
\Crefname{theorem}{Theorem}{Theorems}
\crefname{proposition}{proposition}{propositions}
\Crefname{proposition}{Proposition}{Propositions}
\crefname{corollary}{corollary}{corollaries}
\Crefname{corollary}{Corollary}{Corollaries}
\crefname{lemma}{lemma}{lemmas}
\Crefname{lemma}{Lemma}{Lemmas}
\crefname{definition}{definition}{definitions}
\Crefname{definition}{Definition}{Definitions}
\crefname{conjecture}{conjecture}{conjectures}
\Crefname{conjecture}{Conjecture}{Conjectures}
\crefname{example}{example}{examples}
\Crefname{example}{Example}{Examples}
\crefname{remark}{remark}{remarks}
\Crefname{remark}{Remark}{Remarks}

\newcommand{\Z}{\mathbb{Z}}
\newcommand{\N}{\mathbb{N}}

\newcommand{\C}{\mathbb{C}}

\newcommand{\innprod}[2]{\left\langle \thinspace #1 \thinspace \middle| \thinspace #2 \thinspace \right\rangle}

\newcommand{\conj}[1]{\overline{#1}} 

\newcommand{\norm}[1]{ \left\| #1 \right\| }

\newcommand{\setbuilder}[2] { \left\{ #1 \enskip \middle| \enskip #2 \right\} }

\newcommand{\abs}[1]{ \left| #1 \right| }

\newcommand{\card}[1]{\left| #1 \right|}

\newcommand{\closure}[1]{\overline{#1}}

\newcommand{\set}[1]{\left\{#1\right\}}

\newcommand{\directsum}{\oplus}
\newcommand{\isoto}{\cong}
\newcommand{\boundary}{\partial}
\newcommand{\defeq}{\coloneqq}

\newcommand{\wtilde}[1]{\widetilde{#1}}

\newcommand{\actson}{\curvearrowright}
\newcommand{\injectsinto}{\hookrightarrow}

\newcommand{\tensor}{\otimes}

\DeclareMathOperator{\conv}{conv}

\DeclareMathOperator{\ran}{ran}

\DeclareMathOperator{\supp}{supp}

\DeclareMathOperator{\id}{id}

\DeclareMathOperator{\Fix}{Fix}
\newcommand{\normal}{\triangleleft}
\newcommand{\F}{\mathbb{F}}

\DeclareMathOperator{\Char}{Char}
\DeclareMathOperator{\Ad}{Ad}
\DeclareMathOperator{\Aut}{Aut}

\newcommand{\cF}{\mathcal{F}}

\title[Traces on noncommutative crossed products]{Characterizing traces on crossed products of noncommutative C*-algebras}
\author{Dan Ursu}
\address{Department of Pure Mathematics\\University of Waterloo\\
	200 University Avenue West, Waterloo, Ontario, N2L 3G1, Canada}
\email{dursu@uwaterloo.ca}
\thanks{This work was supported by the Natural Sciences and Engineering Research Council of Canada (NSERC) [grant number PGSD3-535032-2019]. Cette recherche a \'{e}t\'{e} financ\'{e}e par le Conseil de recherches en sciences naturelles et en g\'{e}nie du Canada (CRSNG) [num\'{e}ro de subvention PGSD3-535032-2019].}

\begin{document}
	
	\begin{abstract}
		We give complete descriptions of the tracial states on both the universal and reduced crossed products of a C*-dynamical system consisting of a unital C*-algebra and a discrete group. In particular, we also answer the question of when the tracial states are in canonical bijection with the invariant tracial states on the original C*-algebra. This generalizes the unique trace property for discrete groups. The analysis simplifies greatly in various cases, for example when the conjugacy classes of the original group are all finite, and in other cases gives previously known results, for example when the original C*-algebra is commutative. We also obtain results and examples in the case of abelian groups that contradict existing results in the literature of B{\'e}dos and Thomsen. Specifically, we give a finite-dimensional counterexample, and provide a correction to the result of Thomsen.
	\end{abstract}

	\maketitle
	
	\blfootnote{\copyright~2021. This manuscript version is made available under the CC-BY-NC-ND 4.0 license \\ \url{https://creativecommons.org/licenses/by-nc-nd/4.0/}}

	\tableofcontents

	\section{Introduction and statement of main results}
	
	Understanding the tracial states on a C*-algebra is often of interest, for example in classification theory. In this paper, we concern ourselves with both the reduced and universal crossed products arising from C*-dynamical systems consisting of a unital C*-algebra $A$ and a discrete group $G$. Given any $G$-invariant tracial state on $A$, we give complete descriptions of the tracial extensions to the crossed products. We also translate our characterization into an equivalent condition for when the tracial extension is unique. Finally, in various special cases, we simplify this condition on uniqueness of tracial extension.
	
	To establish notation, $A$ will always denote a unital C*-algebra, and $G$ a discrete group acting on $A$ by *-automorphisms. In addition, the term ``automorphism'' will always mean *-automorphism. The reduced crossed product of this action will be denoted by $A \rtimes_\lambda G$, and the universal crossed product by $A \rtimes_u G$. The unitary corresponding to $t \in G$ in $A \rtimes_\lambda G$ will be denoted by $\lambda_t$, and in $A \rtimes_u G$ will simply be denoted by $t$. Furthermore, $T(A)$ will denote the set of all tracial states on $A$, and $T_G(A)$ the set of tracial states that are invariant under the action of $G$. Finally, ``trace'' will only be used to refer to tracial states.
	
	A key idea in our paper takes inspiration from one of the techniques used by Kennedy and Schafhauser in \cite{kennedy_schafhauser_noncommutative_crossed_products}. In their paper, they study the intersection property of reduced crossed products, i.e.\ the property that every nonzero ideal of $A \rtimes_\lambda G$ has nonzero intersection with $A$. A key point in their paper is that what they call \textit{pseudoexpectations}, introduced in \cite[Section~6]{kennedy_schafhauser_noncommutative_crossed_products}, can be used to characterize the intersection property. These are $G$-equivariant, unital, completely positive maps $\phi : A \rtimes_\lambda G \to I_G(A)$ with $\phi|_A$ being the identity map, where $I_G(A)$ is the $G$-injective envelope of $A$. It is worth noting that this notion of pseudoexpectation is based on the original, different notion of pseudoexpectation introduced by Pitts in \cite{pitts_structure_for_regular_inclusions_I},
	and studied by both Pitts and Zarikian in subsequent papers.
	
	We adapt the notion of pseudoexpectations used by Kennedy and Schafhauser to one that can instead be used to characterize tracial extensions of $\tau \in T_G(A)$ to both the universal and reduced crossed products. It has been previously recognized (for example in \cite[Section~2]{bedos_on_the_uniqueness_of_trace}) that the dynamics of $G$ on $\pi(A)''$, where $\pi : A \to B(H_\tau)$ is the GNS representation of $\tau$, plays an important role in determining the tracial extensions to the crossed products. See the review we give in \Cref{subsection_tracial_gns} for why we have an action of $G$ on this von Neumann algebra, along with other basic properties. This seems to suggest that $\pi(A)''$ is the appropriate object to consider in place of the $G$-injective envelope $I_G(A)$.
	
	As we will make use of the amenable radical of $G$ when dealing with the reduced crossed product, the notion of pseudoexpectation that we introduce works relative to any normal subgroup of $G$. Note that, given a normal subgroup $N \normal G$, we canonically have an action of $G$ on $A \rtimes_u N$ satisfying $s \cdot (at) = (s \cdot a) (sts^{-1})$ for $s \in G$ and $t \in N$, by the universal property of $A \rtimes_u N$.
	
	We also note ahead of time that we will not use the term ``pseudoexpectation'', which as mentioned above refers to maps whose codomain is some kind of injective envelope. Instead, for our purposes, we will adopt the term ``weak expectation'', which instead typically refers to maps whose codomain is an enveloping von Neumann algebra. For example, the usual notion of a weak expectation of a unital inclusion of C*-algebras $A \subseteq B$ is a unital and completely positive map $F : B \to A^{**}$ extending the canonical inclusion $\iota : A \injectsinto A^{**}$. In our case, we replace the universal enveloping von Neumann algebra $A^{**}$ with the von Neumann algebra $\pi(A)''$ generated under the GNS representation $\pi : A \to B(H_\tau)$ of our fixed trace $\tau \in T_G(A)$. Note that the canonical map $\pi : A \to \pi(A)''$ is not necessarily faithful, contrary to the case of $A^{**}$, $I(A)$, or $I_G(A)$.
	
	\begin{definition}
		Let $\tau \in T_G(A)$, let $\pi : A \to B(H_\tau)$ be the GNS representation of $(A,\tau)$, let $M = \pi(A)''$, and let $N \normal G$ be a normal subgroup. A map $F : A \rtimes_u N \to M$ is called a \textit{weak expectation} for $(A,\tau,G,N)$ if it is unital, completely positive, $G$-equivariant, and satisfies $F|_A = \pi$. If $N = G$, then we call such a map a weak expectation for $(A,\tau,G)$.
	\end{definition}
	
	\begin{theorem}
	\label{crossed_product_traces_universal}
		Let $\tau \in T_G(A)$, let $\pi : A \to B(H_\tau)$ be the GNS representation of $(A,\tau)$, let $M = \pi(A)''$, and let $\tau_M$ denote the corresponding faithful normal trace on $M$. Then the following sets are in natural bijection with each other:
		\begin{enumerate}
			\item The set of all weak expectations $F : A \rtimes_u G \to M$ for $(A,\tau,G)$.
			\item The set of all $\set{x_t}_{t \in G} \subseteq M$ satisfying:
			\begin{enumerate}
				\item $x_e = 1$.
				\item $x_t y = (t \cdot y) x_t$ for all $y \in M$ and $t \in G$.
				\item $s \cdot x_t = x_{sts^{-1}}$ for all $s,t \in G$.
				\item The matrix $[x_{st^{-1}}]_{s,t \in \cF}$ is positive for all finite $\cF \subseteq G$.
			\end{enumerate}
			\item $\setbuilder{\sigma \in T(A \rtimes_u G)}{\sigma|_A = \tau}$.
		\end{enumerate}
		The natural map from (1) to (2) is given by letting $x_t = F(t)$, and the natural map from (2) to (3) is given by defining a trace $\sigma \in T(A \rtimes_u G)$ by $\sigma(at) = \tau_M(\pi(a)x_t)$.
	\end{theorem}

	For the case of the reduced crossed product, we replace almost all instances of $G$ with the amenable radical $R_a(G)$, which is the largest amenable normal subgroup of $G$. This was originally introduced by Day in \cite[Section~4, Lemma~1]{day_amenable_semigroups}. The main idea making the case of the reduced crossed product tractable is that, by a result of Bryder and Kennedy \cite[Theorem~5.2]{bryder_kennedy_twisted_crossed_products}, any trace on $A \rtimes_\lambda G$ concentrates on $A \rtimes_\lambda R_a(G)$, in the sense that it vanishes on $a\lambda_t$ whenever $t \notin R_a(G)$. However, $A \rtimes_\lambda R_a(G) = A \rtimes_u R_a(G)$ by amenability of $R_a(G)$, and so we may apply the results we obtained in the case of universal crossed products.

	\begin{theorem}
	\label{crossed_product_traces_reduced}
		Let $\tau \in T_G(A)$, let $\pi : A \to B(H_\tau)$ be the GNS representation of $(A,\tau)$, let $M = \pi(A)''$, and let $\tau_M$ denote the corresponding faithful normal trace on $M$. Then the following sets are in natural bijection with each other:
		\begin{enumerate}
			\item The set of all weak expectations $F : A \rtimes_u R_a(G) \to M$ for $(A,\tau,G,R_a(G))$.
			\item The set of all $\set{x_t}_{t \in R_a(G)} \subseteq M$ satisfying:
			\begin{enumerate}
				\item $x_e = 1$.
				\item $x_t y = (t \cdot y) x_t$ for all $y \in M$ and $t \in R_a(G)$.
				\item $s \cdot x_t = x_{sts^{-1}}$ for all $s \in G$ and $t \in R_a(G)$.
				\item The matrix $[x_{st^{-1}}]_{s,t \in \cF}$ is positive for all finite $\cF \subseteq R_a(G)$.
			\end{enumerate}
			\item $\setbuilder{\sigma \in T(A \rtimes_\lambda G)}{\sigma|_A = \tau}$.
		\end{enumerate}
		The natural map from (1) to (2) is given by letting $x_t = F(\lambda_t)$, and the natural map from (2) to (3) is given by defining a trace $\sigma \in T(A \rtimes_\lambda G)$ by $\sigma(a\lambda_t) = \tau_M(\pi(a)x_t)$ for $t \in R_a(G)$, and $\sigma(a \lambda_t) = 0$ for $t \notin R_a(G)$.
	\end{theorem}
	
	Traces on $A \rtimes_u G$ and $A \rtimes_\lambda G$ are easiest to understand when they correspond exactly to $G$-invariant traces on $A$. Let $E : A \rtimes_\lambda G \to A$ denote the canonical expectation. There is also a canonical expectation from $A \rtimes_u G$ to $A$ given by composing $E$ with the canonical *-homomorphism from $A \rtimes_u G$ to $A \rtimes_\lambda G$.
		
	\begin{remark}
		Any trace $\sigma \in T(A \rtimes_\lambda G)$ satisfies $\sigma|_A \in T_G(A)$. Conversely, any $\tau \in T_G(A)$ gives rise to a trace on $A \rtimes_\lambda G$ by composing with the canonical expectation. That is, $\tau \circ E \in T(A \rtimes_\lambda G)$. Analogous results hold for the universal crossed product.
	\end{remark}
	
	Keeping the above in mind, we generalize the notion of the unique trace property for discrete groups. Recall that $G$ is said to have the unique trace property if the only trace on the reduced group C*-algebra $C^*_\lambda(G) \subseteq B(\ell^2(G))$ is the canonical one, given by $\tau_\lambda(a) = \innprod{a \delta_e}{\delta_e}$. This was shown in \cite[Corollary~4.3]{breuillard_kalantar_kennedy_ozawa_c_simplicity} to be equivalent to $G$ having trivial amenable radical.
	
	\begin{definition}
		Given any $\sigma \in T(A \rtimes_\lambda G)$, we will say that $\sigma$ is \textit{canonical} if it is of the form $\sigma = \tau \circ E$ for some $\tau \in T_G(A)$, or equivalently, if $\sigma = \sigma \circ E$. Canonical traces on $A \rtimes_u G$ are defined analogously. Given $\tau \in T_G(A)$, we will say that it has \textit{unique tracial extension to $A \rtimes_\lambda G$} (or $A \rtimes_u G$) if the only $\sigma \in T(A \rtimes_\lambda G)$ (respectively, $T(A \rtimes_u G)$) satisfying $\sigma|_A = \tau$ is the canonical one.
	\end{definition}

	Before proceeding further, we note that setting $A = \C$ in both of the above theorems indeed gives back previously known results.

	\begin{remark}
		 Setting $A = \C$ in \Cref{crossed_product_traces_universal} gives back the well-known result that traces on the universal group C*-algebra $C^*_u(G)$ correspond to positive definite functions $f : G \to \C$ that are constant on conjugacy classes and satisfy $f(e) = 1$. In addition, setting $A = \C$ in \Cref{crossed_product_traces_reduced} gives back the fact that the unique trace property for the reduced group C*-algebra $C^*_\lambda(G)$ is equivalent to having $R_a(G) = \set{e}$, as in set (2), we may always let $x_t = 1$ for $t \in R_a(G)$.
	\end{remark}

	In the above theorems in set (2), condition (b) in particular highlights a link with proper outerness of the action (see the review in \Cref{subsection_properly_outer_automorphisms} for definitions). Thus, we obtain an immediate corollary which is perhaps a slight generalization of some previously known results (see \cite[Proposition~9]{bedos_simple_unique_trace}, for example).
	
	\begin{corollary}
	\label{properly_outer_implies_unique_trace_property}
		Let $\tau \in T_G(A)$, let $\pi : A \to B(H_\tau)$ be the GNS representation of $(A,\tau)$, and let $M = \pi(A)''$. If the action of $R_a(G)$ on $M$ is properly outer, then $\tau$ has unique tracial extension to $A \rtimes_\lambda G$. If the action of $G$ on $M$ is properly outer, then $\tau$ has unique tracial extension to $A \rtimes_u G$.
	\end{corollary}
	
	Another key idea in our paper is developed in \Cref{section_partial_inner_actions_and_vanishing_obstruction}, where we convert the conditions in set (2) in \Cref{crossed_product_traces_universal,crossed_product_traces_reduced} into conditions on what we call a \textit{partial almost inner action}, with an optional property which we call \textit{positively compatible}. This is similar to the notion of a partial group representation---see, for example, the book of Exel \cite{exel_partial_dynamical_systems}---together with its applications in the work done by Kennedy and Schafhauser in \cite{kennedy_schafhauser_noncommutative_crossed_products}.
	We adapt the notion of a properly outer action (for single automorphisms), giving us what we call a \textit{jointly almost properly outer action}. This condition is again sufficient, and in some special cases necessary, for $\tau \in T_G(A)$ to have unique tracial extension. All of the definitions required for this theorem and its corollaries can be found in \Cref{definition_partial_inner_action,definition_positively_compatible}.
	
	\begin{theorem}
	\label{unique_trace_property_almost_inner_actions}
		Let $\tau \in T_G(A)$, let $\pi : A \to B(H_\tau)$ be the GNS representation of $(A,\tau)$, and let $M = \pi(A)''$.
		\begin{enumerate}
			\item The trace $\tau$ has unique tracial extension to the reduced crossed product $A \rtimes_\lambda G$ if and only if the action of $G$ on $M$ is not partially almost inner relative to the normal subgroup $R_a(G)$ with respect to some nontrivial positively compatible $\set{(p_t, u_t)}_{t \in R_a(G)}$. In particular, it is sufficient for the action to be jointly almost properly outer relative to $R_a(G)$.
			
			\item The trace $\tau$ has unique tracial extension to the universal crossed product $A \rtimes_u G$ if and only if the action of $G$ on $M$ is not partially almost inner with respect to some nontrivial positively compatible $\set{(p_t, u_t)}_{t \in G}$. In particular, it is sufficient for the action to be jointly almost properly outer.
		\end{enumerate}		
	\end{theorem}
	
	The rest of our results are simplifications of the above theorem in certain special cases. Recall that the FC center of a group $G$ is the set of all elements of $G$ with finite conjugacy classes. An FC group is a group in which every conjugacy class is finite, i.e.\ one that is equal to its FC center. It is known that FC groups are amenable,
	and so in particular this next result applies to such groups.
		
	\begin{corollary}
	\label{unique_trace_property_FC_groups}
		Assume $G$ is a group with the property that the amenable radical $R_a(G)$ and the FC center coincide. Let $\tau \in T_G(A)$, let $\pi : A \to B(H_\tau)$ denote the GNS representation, and let $M = \pi(A)''$. Then $\tau$ has unique tracial extension to $A \rtimes_\lambda G$ if and only if the action of $G$ on $M$ is jointly almost properly outer relative to the normal subgroup $R_a(G)$.
	\end{corollary}

	The conditions of the above theorem simplify even further in the case of groups whose amenable radical is equal to the center. In particular, the following corollary applies to abelian groups.

	\begin{corollary}
	\label{unique_trace_property_abelian_groups}
		Assume $G$ is a group with the property that $R_a(G) = Z(G)$. Let $\tau \in T_G(A)$, let $\pi : A \to B(H_\tau)$ denote the GNS representation, and let $M = \pi(A)''$. Then $\tau$ has unique tracial extension to $A \rtimes_\lambda G$ if and only if for any $t \in R_a(G) \setminus \set{e}$, there does not exist a central projection $p \neq 0$ in $M$ and $u \in U(Mp)$ with the properties that:
		\begin{enumerate}
			\item $s \cdot p = p$ and $s \cdot u = u$ for all $s \in G$.
			\item $t$ acts by $\Ad u$ on $Mp$.
		\end{enumerate}
	\end{corollary}
	
	In the case of abelian groups, various results of B{\'e}dos and Thomsen on finite factoriality of von Neumann crossed products \cite[Proposition~11]{bedos_simple_unique_trace}, and unique tracial extension for C*-crossed products \cite[Theorem~4.3]{thomsen_trace_bijection}, respectively, already exist in the literature, but they are incorrect. This is investigated in \Cref{subsection_counterexample_finite}, where a finite-dimensional counterexample is given. \Cref{unique_trace_property_abelian_groups} serves as a correction to the result of Thomsen. Interestingly enough, even though \Cref{subsection_counterexample_finite_cyclic_group} gives a counterexample in the case of finite cyclic groups, Thomsen's result still 
	holds in the case of integer actions.
	
	\begin{theorem}
	\label{unique_trace_property_integers}
		Assume $\alpha \in \Aut(A)$, and consider the corresponding action of $\Z$ on $A$. Let $\tau \in T_{\Z}(A)$, let $\pi : A \to B(H_\tau)$ denote the GNS representation, and let $M = \pi(A)''$. Then $\tau$ has unique tracial extension to $A \rtimes_\lambda \Z$ if and only if the action of $\Z$ on $M$ is properly outer.
	\end{theorem}
	
	Another case in which the characterization simplifies is in the case of crossed products of commutative C*-algebras. The case of the universal crossed product is already known---see \cite[Theorem~2.7]{kawamura_hideo_tomiyama_state_extensions}. Essential freeness and its relation to proper outerness are reviewed in \Cref{subsection_properly_outer_automorphisms}.
	
	\begin{corollary}
	\label{unique_trace_property_commutative}
		Assume $G$ acts on a compact Hausdorff space $X$ by homeomorphisms, and $\mu$ is a $G$-invariant Radon probability measure on $X$. Then $\mu$ has unique tracial extension to $C(X) \rtimes_\lambda G$ if and only if the action of $R_a(G)$ on $(X,\mu)$ is essentially free. For the universal crossed product $C(X) \rtimes_u G$, $\mu$ has unique tracial extension if and only if the action of $G$ on $(X,\mu)$ is essentially free.
	\end{corollary}

	\addtocontents{toc}{\SkipTocEntry}
	\section*{Acknowledgments}	
	
	The author would like to thank his supervisor, Matthew Kennedy, for giving detailed comments and suggestions throughout this paper being written. In addition, the author would also like to thank Erik B{\'e}dos, Mehrdad Kalantar, Sven Raum, and Klaus Thomsen for taking the time to look through a near-finished draft of this paper and giving helpful feedback. The author would also like to thank Se-Jin (Sam) Kim and Nicholas Manor for discussions on the validity of \Cref{tracial_c_embedding_gives_tracial_gns_vn_embedding} and \Cref{essentially_free_iff_properly_outer}, respectively.

	\section{Preliminaries}
	\label{section_preliminaries}
	
	\subsection{Tracial GNS representations}
	\label{subsection_tracial_gns}
		
	Throughout this paper, we will make heavy use of passing from a tracial C*-algebra to the von Neumann algebra it generates under the GNS representation. Here, we establish the basic facts that we will use. This first proposition is well-known---see, for example, \cite[Chapter~V, Proposition~3.19]{takesaki_theory_of_operator_algebras_I}.
	
	\begin{proposition}
		Assume $A$ is a unital C*-algebra, and $\tau \in T(A)$. Let $\pi : A \to B(H_\tau)$ denote the GNS representation and let $M = \pi(A)''$. Then there is a faithful normal trace $\tau_M$ on $M$ satisfying $\tau_M \circ \pi = \tau$.
	\end{proposition}

	Observe that the above proposition makes no assumptions on $\tau \in T(A)$ being faithful---it is always the case that $\tau_M \in T(M)$ is faithful. In addition, $\tau_M$ is uniquely determined by normality.

	It is also a basic fact of von Neumann algebras that we do not need to worry about normality when dealing with *-isomorphisms. The following can be found, for example, in \cite[Chapter~III, Corollary~3.10]{takesaki_theory_of_operator_algebras_I}.
	
	\begin{proposition}
		Assume $\pi : M \to N$ is a *-isomorphism of von Neumann algebras. Then $\pi$ is automatically normal and has normal inverse.
	\end{proposition}
	
	It is well-known that any trace-preserving group action on a C*-algebra will extend to the GNS von Neumann algebra. First, we note the following result, the proof of which is straightforward and left as an exercise to the reader.
	
	\begin{lemma}
	\label{isomorphism_of_c_alg_gives_isomorphism_of_gns_vn_alg}
		Assume $A$ and $B$ are unital C*-algebras, $\tau \in T(A)$ and $\sigma \in T(B)$, and $\rho : A \to B$ is a *-isomorphism satisfying $\sigma \circ \rho = \tau$. Let $\pi_\tau : A \to B(H_\tau)$ and $\pi_\sigma : B \to B(H_\sigma)$ denote the GNS representations, let $M = \pi_\tau(A)''$ and $N = \pi_\sigma(B)''$, and let $\tau_M \in T(M)$ and $\tau_N \in T(N)$ the corresponding faithful normal traces. There is a unique *-isomorphism $\wtilde{\rho} : M \to N$ that satisfies $\wtilde{\rho} \circ \pi_\tau = \pi_\sigma \circ \rho$. In addition, $\sigma_N \circ \wtilde{\rho} = \tau_M$.
	\end{lemma}

	\begin{proposition}
		Assume $A$ is a unital C*-algebra and $\tau \in T_G(A)$. Let $\pi : A \to B(H_\tau)$ denote the GNS representation, let $M = \pi(A)''$, and let $\tau_M$ denote the corresponding faithful normal trace on $M$. Letting $\alpha_t : A \to A$ denote the action of $t \in G$ on $A$, there are *-automorphisms $\wtilde{\alpha}_t : M \to M$ satisfying $\wtilde{\alpha}_t \circ \pi = \pi \circ \alpha_t$, and each $\wtilde{\alpha}_t$ is uniquely determined. In addition, $t \mapsto \wtilde{\alpha}_t$ defines a valid group action, and with respect to this action, we have that $\tau_M \in T_G(M)$, and $\pi : A \to M$ is $G$-equivariant.
	\end{proposition}

	\begin{proof}
		Existence and uniqueness of $\wtilde{\alpha}_t : M \to M$ satisfying $\wtilde{\alpha}_t \circ \pi = \pi \circ \alpha_t$ immediately follows from \Cref{isomorphism_of_c_alg_gives_isomorphism_of_gns_vn_alg}. It remains to check that $t \mapsto \wtilde{\alpha}_t$ indeed gives a group homomorphism:
		$$ \wtilde{\alpha}_s (\wtilde{\alpha}_t (\pi(a))) = \wtilde{\alpha}_s \pi(\alpha_t(a)) = \pi(\alpha_s(\alpha_t(a))) = \pi(\alpha_{st}(a)). $$
		Our earlier remark on the uniqueness of these *-automorphisms tells us that $\wtilde{\alpha}_s \circ \wtilde{\alpha}_t = \wtilde{\alpha}_{st}$. Finally, the fact that $\tau_M \circ \wtilde{\alpha}_t = \tau_M$ tells us $\tau_M \in T_G(M)$, and the fact that $\wtilde{\alpha}_t \circ \pi = \pi \circ \alpha_t$ tells us $\pi : A \to M$ is $G$-equivariant.
	\end{proof}

	GNS representations also behave nicely with respect to trace-preserving inclusions of C*-algebras. The following lemma is likely already known---we offer a proof here for convenience.
	
	\begin{proposition}
		\label{tracial_c_embedding_gives_tracial_gns_vn_embedding}
		Assume $A \subseteq B$ is a unital embedding of C*-algebras, and $\tau \in T(B)$. Let $\pi : A \to B(L^2(A,\tau))$ and $\sigma : B \to B(L^2(B,\tau))$ be the GNS representations of $(A,\tau|_A)$ and $(B,\tau)$, respectively, let $M = \pi(A)''$ and $N = \sigma(B)''$, and let $\tau_M \in T(M)$ and $\tau_N \in T(N)$ be the corresponding faithful normal traces. Then we have an embedding $\iota : M \to N$ with the properties that $\iota(M)$ is a von Neumann subalgebra of $N$, $\iota : M \to \iota(M)$ is a normal *-isomorphism with normal inverse, $\iota \circ \pi = \sigma|_A$, and $\tau_N \circ \iota = \tau_M$.
	\end{proposition}
	
	\begin{proof}
		Observe that we canonically have $L^2(A,\tau) \subseteq L^2(B,\tau)$, and let $F : B(L^2(B,\tau)) \to B(L^2(A,\tau))$ denote the compression map. Given that $L^2(A,\tau)$ is $\sigma(A)$-invariant, we have that $F(\sigma(a)) = \pi(a)$ for all $a \in A$. By normality, we have $F(\sigma(A)'') \subseteq \pi(A)''$.
		
		We claim that $F|_{\sigma(A)''} : \sigma(A)'' \to \pi(A)''$ is injective. Observe that $\tau_M \circ F$ and $\tau_N$ agree on $\sigma(A)$, and so by normality, on $\sigma(A)''$. But $\tau_N$ is faithful, and this forces $F$ to be faithful on $\sigma(A)''$.
		
		Surjectivity of $F|_{\sigma(A)''} : \sigma(A)'' \to \pi(A)''$ is also easy enough to deduce---the unit ball of $\sigma(A)''$ is weak*-compact, and hence by normality it maps to a weak*-closed subset of $\pi(A)''$. In addition, the unit ball of $\sigma(A)$ maps to a norm-dense subset of the unit ball of $\pi(A)$ (this is true for any quotient map of C*-algebras). These two facts, combined with Kaplansky density, tell us that the image of the unit ball of $\sigma(A)''$ is the entire unit ball of $\pi(A)''$. Linearity takes care of the rest.
		
		In summary, we have shown that $F|_{\sigma(A)''} : \sigma(A)'' \to \pi(A)''$ is a *-isomorphism. We claim that $\iota \defeq (F|_{\sigma(A)''})^{-1} : \pi(A)'' \to \sigma(A)''$ is the embedding we are looking for. By construction, we have $\iota(\pi(a))$ to $\sigma(a)$. From here, we see that
		$$ \tau_N(\iota(\pi(a))) = \tau_N(\sigma(a)) = \tau(a) = \tau_M(\pi(a)), $$
		and so by normality, $\tau_N \circ \iota$ and $\tau_M$ agree on all of $M$.
	\end{proof}

	\subsection{Properly outer automorphisms}
	\label{subsection_properly_outer_automorphisms}
		

	It has long been recognized that proper outerness of an action of $G$ on a C*-algebra or a von Neumann algebra leads to nice structure theory for the corresponding crossed product, particularly in the von Neumann algebra case---see, for example, \cite[Theorem~3.3]{kallman_free_actions}. This is a generalization of essential freeness for measure spaces. To establish notation, given a set $X$ and a map $\alpha : X \to X$, we denote the set of fixed points $\setbuilder{x \in X}{\alpha(x) = x}$ by $\Fix(\alpha)$.
	
	\begin{definition}
	\label{definition_essentially_free}
		Assume $X$ is a compact Hausdorff space, $\alpha : X \to X$ a homeomorphism, and $\mu$ an $\alpha$-invariant Radon probability measure on $X$. We say that $\alpha$ is \textit{essentially free} on $(X,\mu)$ if $\mu(\Fix(\alpha)) = 0$. If $G$ is a group acting on $X$ by $\mu$-invariant homeomorphisms $\alpha_t$, we say that the action is essentially free on $(X,\mu)$ if each $\alpha_t$ is essentially free for $t \in G \setminus \set{e}$.
	\end{definition}
	
	Kallman introduced in \cite[Definition~1.3]{kallman_free_actions} a notion of freely acting automorphisms for general von Neumann algebras, as opposed to just $L^\infty(X,\mu)$:
	
	\begin{definition}
	\label{definition_freely_acting}
		Let $M$ be a von Neumann algebra and $\alpha \in \Aut(M)$. We say that $\alpha$ is \textit{freely acting} if whenever $xy = \alpha(y)x$ for all $y \in M$, we have $x = 0$.
	\end{definition}

	General automorphisms on von Neumann algebras enjoy a very nice decomposition theory into an inner part and a freely acting part---see \cite[Theorem~1.11]{kallman_free_actions}, along with its proof.
	
	\begin{theorem}
	\label{innner_free_decomposition}
		Let $M$ be a von Neumann algebra and $\alpha \in \Aut(M)$. There is a largest $\alpha$-invariant central projection $p \in M$ with the property that $\alpha|_{Mp}$ is inner. In addition, $\alpha|_{M(1-p)}$ is freely acting. Finally, the decomposition $\alpha = \alpha_1 \directsum \alpha_2$ and $M = M_1 \directsum M_2$, with $\alpha_i \in \Aut(M_i)$ and the property that $\alpha_1$ is inner and $\alpha_2$ is freely acting, is unique.
	\end{theorem}

	Proper outerness is an equivalent formulation of freeness (nowadays, the terms are often used interchangeably),
	and is usually defined as follows:
	
	\begin{definition}
	\label{definition_properly_outer}
		Let $M$ be a von Neumann algebra and $\alpha \in \Aut(M)$. We say that $\alpha$ is \textit{properly outer} if there is no nonzero $\alpha$-invariant central projection $p \in M$ with the property that $\alpha|_{Mp}$ is inner. If $G$ is a group acting on $M$ by *-automorphisms $\alpha_t$, we say that the action is properly outer if each $\alpha_t$ is properly outer for $t \in G \setminus \set{e}$.
	\end{definition}

	We will not make use of this following definition, but it is worth noting that we call a group action $\alpha : G \to \Aut(M)$ \textit{inner} if there is a group homomorphism $\beta : G \to U(M)$ with the property that $\alpha(t) = \Ad \beta(t)$. It is important to keep in mind that this is not equivalent to having each $\alpha(t)$ be inner---indeed, it is not hard to check that the example in \Cref{subsection_counterexample_finite} is in fact a finite-dimensional example of this phenomenon. In addition, we call the action \textit{outer} if each $\alpha(t)$, $t \in G \setminus \set{e}$, is outer. Observe that if $M$ is a factor, then outer and properly outer are equivalent.
	
	This next result is well-known, and highlights the fact that proper outerness truly does generalize the notion of essential freeness. The proof is slightly nontrivial and hard to find in the literature, and so we include it here.
	
	\begin{proposition}
	\label{essentially_free_iff_properly_outer}
		Assume $X$ is a compact Hausdorff space, $\alpha : X \to X$ a homeomorphism, and $\mu$ an $\alpha$-invariant Radon probability measure. Then $\alpha$ is essentially free on $(X,\mu)$ if and only if the corresponding automorphism on $L^\infty(X,\mu)$ is properly outer.
	\end{proposition}

	\begin{proof}
		If $\alpha$ were not essentially free on $(X,\mu)$, then $p = 1_{\Fix(\alpha)}$ is a nonzero $\alpha$-invariant central projection in $L^\infty(X,\mu)$, and the action on $(L^\infty(X,\mu))p$ is trivial.
		
		Conversely, assume such a projection $p \in L^\infty(X,\mu)$ exists, and let $E = \supp p$. Replacing $E$ by $\cup_{n \in \Z} \alpha^n(E)$, we may assume without loss of generality that $E$ is $\alpha$-invariant. We claim that $Y \defeq E \setminus \Fix(\alpha)$ is a null set. Given $y \in Y$, we may choose an open neighborhood $U_y$ with the property that $\alpha(U_y) \cap U_y = \emptyset$. Observe that $U_y \cap Y$ is a null set by our assumption that $\alpha$ acts trivially on $L^\infty(E,\mu)$. Now, $\mu$ is inner regular on all sets (it is outer regular, and we may take complements), and so given any $\varepsilon > 0$, we may choose a compact set $K \subseteq Y$ with $\mu(Y \setminus K) < \varepsilon$. By compactness, $K$ admits a finite subcover from $\set{U_y}_{y \in Y}$, and using the fact that every $U_y \cap K$ is a null set, we deduce that $K$ is a null set. Consequently, so is $Y$. Thus, without loss of generality, we have $E \subseteq \Fix(\alpha)$, and so the action of $\alpha$ on $(X,\mu)$ is not essentially free.
	\end{proof}

	Although we will not make use of this fact, it is worth keeping in mind that ``central'' is often omitted from \Cref{innner_free_decomposition} and \Cref{definition_properly_outer}. It is a result of Borchers, \cite[Lemma~5.7]{borchers_inner_automorphisms},
	that if $e$ is any projection in $M$, not necessarily central, and $\alpha \in \Aut(M)$ satisfies $\alpha(e) = e$ and is inner on $eMe$, then it is inner on $Mp$, where $p$ is the central cover of $e$.

	\section{Almost inner actions}
	\label{section_partial_inner_actions_and_vanishing_obstruction}
		
	This section builds on what was reviewed in \Cref{subsection_properly_outer_automorphisms}. As previously mentioned, we aim to convert the conditions in \Cref{crossed_product_traces_universal,crossed_product_traces_reduced}, set (2), into conditions on inner actions on corners of the von Neumann algebra $M$, which are often much easier to check in practice. This will be done by taking the polar decomposition of the elements $x_t$, and this is the motivation behind the definitions that follow. Observe that condition (b) in these theorems is precisely the identity used in the definition of freely acting automorphisms. Condition (d), however, has no obvious nice resulting condition on the unitaries that we obtain, and so we do not include any analogous condition in our definition of \textit{partially almost inner} below. Instead, we include it as a separate property which we call \textit{positively compatible}.
	
	\begin{definition}
	\label{definition_partial_inner_action}
		Assume $M$ is a von Neumann algebra, $N \normal G$ is normal, and $G$ acts on $M$ by *-automorphisms. We say that the action is \textit{partially almost inner relative to $N$} with respect to $\set{(p_t, u_t)}_{t \in N}$ if:
		\begin{enumerate}
			\item Given any $t \in N$, $p_t$ is a central projection in $M$ satisfying $t \cdot p_t = p_t$, $u_t$ is a unitary in $Mp_t$, and moreover, $t$ acts on $Mp_t$ by $\Ad u_t$.
			\item $p_e = 1$ and $u_e = 1$.
			\item $p_t = p_{t^{-1}}$ and $u_t^* = u_{t^{-1}}$ for all $t \in N$.
			\item $s \cdot p_t = p_{sts^{-1}}$ and $s \cdot u_t = u_{sts^{-1}}$ for all $s \in G$ and $t \in N$.
		\end{enumerate}
		If $p_t = 0$ for all $t \in N \setminus \set{e}$, we call $\set{(p_t, u_t)}_{t \in N}$ \textit{trivial}, and \textit{nontrivial} otherwise. If there exists a choice of $\set{(p_t, u_t)}_{t \in N}$ with $p_t = 1$ for all $t \in N$, then we say that the action is \textit{almost inner \textit{relative to $N$}}. If, in addition, $N = G$, then we simply call the action \textit{almost inner}. We say that the action is \textit{jointly almost properly outer relative to $N$} if the only $\set{(p_t, u_t)}_{t \in N}$ with respect to which it is partially almost inner is the trivial one. We will simply call the action \textit{jointly almost properly outer} if it is jointly properly outer relative to $G$.
	\end{definition}

	\begin{remark}
		Consider the above definition in the case of $N = G$. It is worth noting that if $p_t = 1$ for all $t \in G$, then the map $t \mapsto u_t$ is not necessarily an inner action of $G$ on $M$, i.e.\ we are not guaranteed that $u_{st} = u_s u_t$. This is the motivation behind the term \textit{almost}. In addition, the usual definitions of a properly outer (or just outer) action $\alpha : G \to \Aut(M)$ is that each individual $\alpha(t)$, $t \neq e$, is properly outer (respectively, outer). The term \textit{jointly} highlights the fact that we require compatibility conditions between each of the individual $\alpha(t)$.
	\end{remark}

	\begin{definition}
		\label{definition_positively_compatible}
		Let $\set{(p_t, u_t)}_{t \in N}$ be as in \Cref{definition_partial_inner_action}. We say that $\set{(p_t, u_t)}_{t \in N}$ are \textit{positively compatible} if there exist elements $\set{y_t}_{t \in N} \subseteq M$ such that:
		\begin{enumerate}
			\item $y_t \in Z(M)$ and $y_t \geq 0$ for all $t \in N$.
			\item The projection onto $\closure{\ran} y_t$ is $p_t$ for all $t \in N$.
			\item $y_e = 1$.
			\item $y_t = y_{t^{-1}}$ for all $t \in N$.
			\item $s \cdot y_t = y_{sts^{-1}}$ for all $s \in G$ and $t \in N$.
			\item Given any finite $\cF \subseteq N$, the matrix $[u_{st^{-1}} y_{st^{-1}}]_{s,t \in \cF}$ is positive.
		\end{enumerate}
	\end{definition}

	It is worth noting that positive compatibility is not a redundant condition, as the following example shows:
	
	\begin{example}
		It is known that it is possible to construct an infinite group $G$ with only two conjugacy classes, using HNN extensions. A proof can be found in the original paper by Higman, Neumann, and Neumann---see \cite[Theorem~III]{higman_neumann_neumann_hnn_extensions}. Let $A = \C$, let $p_t = 1$ for all $t \in G$, and let $u_t = -1$ for $t \neq e$. It is clear that the trivial action is almost inner with respect to $\set{(p_t,u_t)}_{t \in G}$. However, we claim that $\set{(p_t,u_t)}_{t \in G}$ is not positively compatible. To this end, assume otherwise and let $\set{y_t}_{t \in G}$ be as in \Cref{definition_positively_compatible}, and observe that $y_t$ are all some positive constant $\gamma > 0$ for $t \neq e$. Now letting $\cF \subseteq G$ be any finite subset with $\card{\cF} = n$, we have that the matrix
		$$
		\begin{bmatrix}
			1 & -\gamma & \dots & -\gamma \\
			-\gamma & 1 & \ddots & \vdots \\
			\vdots  & \ddots & \ddots & -\gamma \\
			-\gamma & \dots & -\gamma & 1
		\end{bmatrix}
		$$
		is positive. Letting $[-\gamma]$ denote the $n \times n$ matrix with all entries being $-\gamma$, the above matrix is equal to $(1+\gamma)I + [-\gamma]$, and basic linear algebra tells us that the eigenvalues of this matrix are $(1+\gamma) - n\gamma$ and $1 + \gamma$. In particular, $(1+\gamma) - n\gamma < 0$ if $n$ is sufficiently large, contradicting the positivity of the above matrix.
	\end{example}
	
	\begin{proposition}
		\label{jointly_properly_outer_jointly_freely_acting_equivalence}
		Assume $M$ is a von Neumann algebra, $N \normal G$ is normal, $G$ acts on $M$ by *-automorphisms, and $\set{x_t}_{t \in N} \subseteq M$ is such that:
		\begin{enumerate}
			\item $x_t y = (t \cdot y) x_t$ for all $y \in M$ and $t \in N$.
			\item $x_e = 1$.
			\item $x_t^* = x_{t^{-1}}$ and $t \in N$.
			\item $s \cdot x_t = x_{sts^{-1}}$ for all $s \in G$ and $t \in N$.
		\end{enumerate}
		Given $t \in N$, consider the polar decomposition of $x_t$, i.e. let $u_t$ be the unique partial isometry such that both $x_t = u_t \abs{x_t}$ and $u_t^*u_t$ is the projection onto $\closure{\ran} \abs{x_t}$, and furthermore denote this projection by $p_t$. Then the action is partially almost inner relative to $N$ with respect to $\set{(p_t, u_t)}_{t \in N}$. Moreover, if for every finite $\cF \subseteq N$, we have that the matrix $[x_{st^{-1}}]_{s,t \in \cF}$ is positive, then $\set{(p_t, u_t)}_{t \in N}$ is positively compatible with respect to $\set{\abs{x_t}}_{t \in N}$.
		
		Conversely, if the action is partially almost inner relative to $N$ with respect to $\set{(p_t, u_t)}_{t \in N}$, then $x_t = u_t$ satisfy the above conditions. If, in addition, $\set{(p_t, u_t)}_{t \in N}$ is positively compatible with respect to $\set{y_t}_{t \in N}$, then $x_t = u_t y_t$ satisfy the above conditions, and also satisfy the property that for any finite $\cF \subseteq N$, the matrix $[x_{st^{-1}}]_{s,t \in \cF}$ is positive.
	\end{proposition}
	
	\begin{proof}
		This proof is somewhat similar to the proof of \cite[Theorem~1.1]{kallman_free_actions}. For convenience, we recreate the necessary parts here. Assume $\set{x_t}_{t \in N}$ is such a collection. Observe that for $w \in U(M)$, we have
		$$ w^*x_t^*x_tw = x_t^* (t \cdot w)^* (t \cdot w) x_t = x_t^*x_t, $$
		which shows $x_t^*x_t \in Z(M)$. Thus, we have that $\abs{x_t}$ and $p_t \in W^*(\abs{x_t})$ also lie in the center. Given $x_t$ is fixed by $t$, so is $p_t$. In addition,
		$$ x_t x_t^* = (t \cdot x_t^*) x_t = x_t^* x_t, $$
		i.e.\ $x_t$ is normal. Now, the equality
		$$ \closure{\ran} x_t = (\ker x_t^*)^\perp = (\ker \abs{x_t^*})^\perp = \closure{\ran} \abs{x_t^*} = \closure{\ran} \abs{x_t} $$
		tells us that $u_t u_t^* = u_t^* u_t$, i.e.\ $u_t$ is a unitary in $Mp_t$. Furthermore, we note that
		$$ u_t y \abs{x_t} = u_t \abs{x_t} y = x_t y = (t \cdot y) x_t = (t \cdot y) u_t \abs{x_t}, $$
		which shows that $t$ acts by $\Ad u_t$ on $Mp_t$. This gives us property (1) of partial almost inner actions. Furthermore, it is clear that we have both $u_e = 1$ and $p_e = 1$ (property (2)).
		
		Given that $\abs{x_{t^{-1}}} = \abs{x_t^*} = \abs{x_t}$, we have $p_{t^{-1}} = p_t$. Now observe that
		$$ u_t^* \abs{x_t} = \abs{x_t}^* u_t^* = x_t^* = x_{t^{-1}} = u_{t^{-1}} \abs{x_{t^{-1}}} = u_{t^{-1}} \abs{x_t}. $$
		Given that $u_t^*$ and $u_{t^{-1}}$ share the same initial projection $p_t$, it follows from uniqueness of polar decomposition that $u_t^* = u_{t^{-1}}$. This is property (3).
		
		Finally, given $s \in G$ and $t \in N$, we see that
		$$ u_{sts^{-1}} \abs{x_{sts^{-1}}} = x_{sts^{-1}} = s \cdot x_t = (s \cdot u_t)(s \cdot \abs{x_t}) = (s \cdot u_t) \abs{x_{sts^{-1}}}. $$
		We wish to conclude that $u_{sts^{-1}}$ and $s \cdot u_t$ have the same initial projection. By the crossed product construction, we may assume without loss of generality that $M \subseteq B(H)$, where $G$ acts on $H$ by unitaries $\lambda_t$, and $t \cdot y = \lambda_t y \lambda_t^*$ for all $y \in M$. This gives us that $(s \cdot u_t)^*(s \cdot u_t) = \lambda_s (u_t^*u_t) \lambda_s^*$ is the projection onto
		$$ \lambda_s(\closure{\ran} \abs{x_t}) = \closure{\ran}(\lambda_s \abs{x_t}) = \closure{\ran}(\lambda_s \abs{x_t} \lambda_s^*) = \closure{\ran} \abs{x_{sts^{-1}}}. $$
		Thus, $u_{sts^{-1}}$ and $s \cdot u_t$ share the same initial projection $p_{sts^{-1}}$. Uniqueness of polar decomposition tells us that they are therefore equal, which also gives us that $s \cdot p_t = p_{sts^{-1}}$. This is property (4).
		
		If $[x_{st^{-1}}]_{s,t \in \cF}$ is positive for all finite $\cF \subseteq N$, then it follows immediately from the definition and from the work that was done above that $\set{(p_t, u_t)}_{t \in N}$ are positively compatible with respect to $\set{\abs{x_t}}_{t \in N}$.
		
		The converse given for converting $\set{(p_t, u_t)}_{t \in N}$ back into elements $x_t$ satisfying the given properties follows from the definitions and is straightforward to verify.
	\end{proof}

	The intersection property for noncommutative reduced crossed products is studied in \cite{kennedy_schafhauser_noncommutative_crossed_products}. Their results show that if the action $G \actson I_G(A)$, where $I_G(A)$ denotes the $G$-injective envelope of $A$, is properly outer, then $A \rtimes_\lambda G$ has the intersection property. Moreover, if the action $G \actson I(A)$,
	where $I(A)$ denotes the usual injective envelope of $A$, has a property they call \textit{vanishing obstruction}, then the converse to this result holds. Here, we show that a very mild adaptation of the intersection property is enough to guarantee that a partial almost inner action is positively compatible.
	
	\begin{proposition}
	\label{vanishing_obstruction_implies_positively_compatible}
		Assume $M$ is a von Neumann algebra, $N \normal G$ is normal, $G$ acts on $M$ by *-automorphisms, and the action is partially almost inner relative to $N$ with respect to $\set{(p_t, u_t)}_{t \in N}$. If, in addition, we have that $p_s p_t \leq p_{st}$ and $u_s u_t = u_{st} p_s p_t$ for all $s,t \in N$, then $\set{(p_t, u_t)}_{t \in N}$ is positively compatible with respect to $\set{p_t}_{t \in N}$.
	\end{proposition}

	\begin{proof}
		Let $\set{s_1, \dots, s_n}$ be a finite subset of $N$. We wish to show that the matrix $[u_{s_i s_j^{-1}}]$ is positive. Let $Z(M) = C(X)$, and consider the sets $\supp p_{s_i s_j^{-1}} \subseteq X$. We may choose finitely many disjoint sets $E_k \subseteq X$ such that $\sqcup_k E_k = X$, and for any $i$, $j$, and $k$, we have $E_k \subseteq \supp p_{s_i s_j^{-1}}$ or $E_k \cap \supp p_{s_i s_j^{-1}} = \emptyset$. These sets can be chosen to be finite intersections of sets of the form $\supp p_{s_i s_j^{-1}}$ and their complements, making each $E_k$ clopen. We will prove that $[u_{s_i s_j^{-1}} 1_{E_k}] \geq 0$ for every $k$.
		
		To this end, fix $k$, and define a relation on $\set{1, \dots, n}$ by $i \sim j$ if and only if $E_k \subseteq \supp p_{s_i s_j^{-1}}$. This is in fact an equivalence relation---it is clear that this is reflexive and symmetric. Transitivity follows from the fact that $p_{s_{i_1} s_{i_2}^{-1}} p_{s_{i_2} s_{i_3}^{-1}} \leq p_{s_{i_1} s_{i_3}^{-1}}$. Thus, if we assume without loss of generality that $\set{s_1, \dots, s_n}$ are ordered such that the equivalence classes are of the form $\set{m, m+1, \dots, m+l}$, then $[u_{s_i s_j^{-1}} 1_{E_k}]$ becomes a block diagonal matrix, where each block of the diagonal is of the form $[u_{s_i s_j^{-1}} 1_{E_k}]_{i,j=m,\dots,m+l}$, and $E_k \subseteq \supp p_{s_i s_j^{-1}}$ for every element in this submatrix. Hence, to prove our original matrix is positive, we may assume without loss of generality that $E_k \subseteq \supp p_{s_i s_j^{-1}}$ for all $i$ and $j$. This matrix is positive, as
		$$
		\begin{bmatrix}
			u_{s_1 s_1^{-1}} 1_{E_k} \\
			\vdots \\
			u_{s_n s_1^{-1}} 1_{E_k}
		\end{bmatrix}
		\begin{bmatrix}
			u_{s_1 s_1^{-1}} 1_{E_k} \\
			\vdots \\
			u_{s_n s_1^{-1}} 1_{E_k}
		\end{bmatrix}^*
		=
		[u_{s_i s_j^{-1}} 1_{E_k}].
		$$
	\end{proof}

	\section{Proof of main results}
	
	As before, $A$ denotes a unital C*-algebra and $G$ a discrete group acting on $A$ by *-automorphisms. Throughout this section, we will fix an invariant trace $\tau \in T_G(A)$, denote by $\pi : A \to B(H_\tau)$ the GNS representation, let $M = \pi(A)''$, and let $\tau_M$ be the corresponding faithful normal trace on $M$.
	
	This first lemma is likely already known, and we give a quick proof for convenience.
	We will denote $A[G] \defeq \set{\sum_{\text{finite}} a_t w_t}$, i.e.\ the set of finitely-supported functions from $G$ to $A$, together with the usual *-algebraic operations obtained by viewing this as a subset of $A \rtimes_u G$. A function $\phi : A[G] \to \C$ is said to be \textit{positive definite} if for any $f \in A[G]$, we have $\phi(f^*f) \geq 0$.
	
	\begin{lemma}
	\label{positive_definite_function_on_AG_extends_to_state}
		Assume $\phi : A[G] \to \C$ is a positive definite function satisfying $\phi(1) = 1$. Then $\phi$ extends to a state on $A \rtimes_u G$.
	\end{lemma}

	\begin{proof}
		This proof is essentially a modified GNS construction. Define a sesquilinear form on $A[G]$ by $\innprod{f_1}{f_2} \defeq \phi(f_2^*f_1)$, and observe that this is positive, as $\phi$ is positive definite. Letting $N = \setbuilder{f \in A[G]}{\innprod{f}{f} = 0}$, we have that the completion of $A[G]/N$ with respect to the corresponding quotient inner product becomes a Hilbert space, which we will denote by $H$.
		
		It is clear that we have a unitary representation $u : G \to U(H)$ given by $u(s) f = w_s f$ for $f \in A[G]$. We also have a *-representation $\rho : A \to B(H)$ given by $\rho(a) f = af$, as
		$$ \innprod{af}{af} = \phi(f^*a^*af) \leq \norm{a}^2 \phi(f^*f), $$
		where this last equality holds due to the fact that
		$$ \phi(f^*(\norm{a}^2 - a^*a)f) = \phi(f^*(\norm{a}^2 - a^*a)^{1/2} (\norm{a}^2 - a^*a)^{1/2} f) \geq 0. $$
		(It is a subtle but important point that $\norm{a}^2 - a^*a$ still admits a positive square root in $A[G]$. This is not necessarily true anymore if we replace $a$ with an arbitrary element of $A[G]$).
		Moreover, $\rho$ and $u$ form a covariant pair. By the universal property of $A \rtimes_u G$, we obtain a *-homomorphism $\wtilde{\rho} : A \rtimes_u G \to B(H)$ given by $\wtilde{\rho}(as) = \rho(a) u(s)$. Consequently, we obtain a positive functional $\sigma \in (A \rtimes_u G)^*$ given by $\sigma(at) = \innprod{\wtilde{\rho}(at) w_e}{w_e} = \phi(aw_t)$. It is clear that, in addition, $\sigma(1) = 1$.
	\end{proof}

	\begin{proposition}
	\label{x_t_trace_existence}
		Assume $\set{x_t}_{t \in G} \subseteq M$ satisfy the assumptions of \Cref{crossed_product_traces_universal}, set (2). Then there is a trace $\sigma \in T(A \rtimes_u G)$ satisfying $\sigma(at) = \tau_M(\pi(a) x_t)$.
	\end{proposition}

	\begin{proof}
		In light of \Cref{positive_definite_function_on_AG_extends_to_state}, to show that we at least obtain a state $\sigma \in S(A \rtimes_u G)$ with the above property, it suffices to show that the function $\sigma : A[G] \to \C$ given by $\sigma(aw_t) = \tau_M(\pi(a) x_t)$ is positive definite. To this end, assume $f = \sum_{i=1}^n a_{s_i} w_{s_i} \in A[G]$. We have that
		\begin{align*}
			\sigma(f^*f) &= \sum_{i=1}^n \sum_{j=1}^n \sigma(w_{s_i}^* a_{s_i}^* a_{s_j} w_{s_j}) \\
			&= \sum_{i=1}^n \sum_{j=1}^n \sigma((s_i^{-1} \cdot (a_{s_i}^* a_{s_j})) w_{s_i^{-1} s_j}) \\
			&= \sum_{i=1}^n \sum_{j=1}^n \tau_M((s_i^{-1} \cdot \pi(a_{s_i})^*) (s_i^{-1} \cdot \pi(a_{s_j})) x_{s_i^{-1} s_j}) \\
			&= \tau_M\left(\sum_{i=1}^n \sum_{j=1}^n (s_i^{-1} \cdot \pi(a_{s_i})^*) x_{s_i^{-1} s_j} (s_j^{-1} \cdot \pi(a_{s_j}))\right).
		\end{align*}
		Observe, however, that
		\begin{align*}
		&\begin{bmatrix}
			s_1^{-1} \cdot \pi(a_{s_1}) \\
			\vdots \\
			s_n^{-1} \cdot \pi(a_{s_n})
		\end{bmatrix}^*
		\begin{bmatrix}
			x_{s_1^{-1} s_1} & \dots & x_{s_1^{-1} s_n} \\
			\vdots & & \vdots \\
			x_{s_n^{-1} s_1} & \dots & x_{s_n^{-1} s_n}
		\end{bmatrix}
		\begin{bmatrix}
			s_1^{-1} \cdot \pi(a_{s_1}) \\
			\vdots \\
			s_n^{-1} \cdot \pi(a_{s_n})
		\end{bmatrix} \\
		&=
		\sum_{i=1}^n \sum_{j=1}^n (s_i^{-1} \cdot \pi(a_{s_i})^*) x_{s_i^{-1} s_j} (s_j^{-1} \cdot \pi(a_{s_j})),
		\end{align*}
		guaranteeing that $\sigma$ is positive definite. It remains to show that the extension to $A \rtimes_u G$ is still a trace:
		\begin{align*}
			\sigma((as) (bt)) &= \sigma((a (s \cdot b)) st) \\
			&= \tau_M(\pi(a) (s \cdot \pi(b)) x_{st}) \\
			&= \tau_M(s^{-1} \cdot (\pi(a) (s \cdot \pi(b)) x_{st})) \\
			&= \tau_M((s^{-1} \cdot \pi(a)) \pi(b) x_{ts}) \\
			&= \tau_M(\pi(b) x_{ts} (s^{-1} \cdot \pi(a))) \\
			&= \tau_M(\pi(b) (t \cdot \pi(a)) x_{ts}) \\
			&= \sigma((b (t \cdot a)) ts) \\
			&= \sigma((bt) (as))
		\end{align*}
	\end{proof}

	\begin{remark}
		The first half of the proof of \Cref{x_t_trace_existence} does not use the fact that $\tau$ (and hence $\tau_M$) is a trace. Thus, if we assume that $\tau$ is only a $G$-invariant state, we still obtain a state $\sigma \in S(A \rtimes_u G)$ given by $\sigma(at) = \tau_M(\pi(a)x_t)$, except $\sigma$ is of course not necessarily a trace anymore.
	\end{remark}

	\begin{lemma}
	\label{x_t_trace_uniqueness}
		Assume $\sigma_1, \sigma_2 \in T(A \rtimes_u G)$ are two states satisfying $\sigma_1(at) = \tau_M(\pi(a) x_t)$ and $\sigma_2(at) = \tau_M(\pi(a) y_t)$ for some $\set{x_t}_{t \in G}, \set{y_t}_{t \in G} \subseteq M$. If $\sigma_1 = \sigma_2$, then $x_t = y_t$ for all $t \in G$.
	\end{lemma}

	\begin{proof}
		Assume otherwise, and fix some $t \in G$ with $x_t \neq y_t$. Letting $(a_\lambda) \subseteq A$ be a net with the property that $(\pi(a_\lambda))$ is weak*-convergent to $(x_t - y_t)^*$, we see that
		$$ (\sigma_1 - \sigma_2)(a_\lambda t) = \tau_M(\pi(a_\lambda)(x_t - y_t)) \to \tau_M((x_t - y_t)^*(x_t - y_t)). $$
		This limit value is nonzero, as $\tau_M$ is faithful. Thus, there is some $\lambda$ such that $\sigma_1(a_\lambda t)$ and $\sigma_2(a_\lambda t)$ differ.
	\end{proof}

	\begin{proof}[Proof of \Cref{crossed_product_traces_universal}]
		Starting with any weak expectation $F : A \rtimes_u G \to M$ for $(A,\tau,G)$ and letting $x_t = F(t)$,
		we note that that $A$ lies in the multiplicative domain of $F$---see, for example, \cite[Proposition~1.5.7]{brown_ozawa_c_algebras_finite_dimensional_approximations}, for a review of multiplicative domain. Thus, $F(at) = \pi(a) x_t$, and so the map between sets (1) and (2) is necessarily injective. It remains to show that $x_t$ indeed satisfy all of the aforementioned properties. We have that $x_e = 1$ follows from $F$ being unital, and $s \cdot x_t = x_{sts^{-1}}$ follows from $F$ being $G$-equivariant. Now, given any $a \in A$ and $t \in G$, observe that
		$$ x_t \pi(a) = F(ta) = F((t \cdot a) t) = (t \cdot \pi(a)) x_t. $$
		Given that $\pi(A)$ is weak*-dense in $M$, taking limits allows us to conclude that $x_t y = (t \cdot y) x_t$ holds for all $y \in M$. Finally, given ${s_1, \dots, s_n} \subseteq G$, we note that
		$$
		F^{(n)}\left(
		\begin{bmatrix}
			s_1 \\
			\vdots \\
			s_n
		\end{bmatrix}
		\begin{bmatrix}
			s_1 \\
			\vdots \\
			s_n
		\end{bmatrix}^*
		\right)
		=
		\begin{bmatrix}
			x_{s_1 s_1^{-1}} & \dots & x_{s_1 s_n^{-1}} \\
			\vdots & & \vdots \\ 
			x_{s_n s_1^{-1}} & \dots & x_{s_n s_n^{-1}}
		\end{bmatrix}.
		$$
		Complete positivity of $F$ says that $[x_{s_i s_j^{-1}}]$ is therefore positive.
		
		Now, starting with any $\set{x_t}_{t \in G} \subseteq M$ as in (2), \Cref{x_t_trace_existence} tells us that $\sigma(at) = \tau_M(\pi(a) x_t)$ indeed defines a valid trace. Moreover, this map from (2) to (3) is injective by \Cref{x_t_trace_uniqueness}.
		
		Finally, we can show the maps from sets (1) to (2) and (2) to (3) are bijective by showing that their composition is surjective. That is, we need to show that for any $\sigma \in T(A \rtimes_u G)$ satisfying $\sigma|_A = \tau$, there exist some weak expectation $F : A \rtimes_u G \to M$ for $(A,\tau,G)$ such that $\sigma = \tau_M \circ F$.
		
		To this end, fix such a $\sigma$, let $\rho : A \rtimes_u G \to B(H_\rho)$ be the GNS representation of $(A \rtimes_u G, \sigma)$, let $N = \rho(A \rtimes_u G)''$, and let $\sigma_N$ denote the corresponding faithful normal trace on $N$. Given that $(A,\tau) \subseteq (A \rtimes_u G, \sigma)$ is a trace-preserving embedding, this canonically gives a trace-preserving embedding $(M,\tau_M) \subseteq (N,\sigma_N)$ sending $\pi(a)$ to $\sigma(a)$ by \Cref{tracial_c_embedding_gives_tracial_gns_vn_embedding}. There is a unique normal conditional expectation $F' : N \to M$ satisfying $\sigma_N = \tau_M \circ F'$---see, for example, \cite[Lemma~1.5.11]{brown_ozawa_c_algebras_finite_dimensional_approximations}. We let $F = F' \circ \rho$, and show that this is the map we are looking for. Observe that
		$$ \tau_M(F(at)) = \tau_M(F'(\rho(at))) = \sigma_N(\rho(at)) = \sigma(at), $$
		i.e.\ $\tau_M \circ F = \sigma$. The only non-trivial fact remaining is to show that $F$ is $G$-equivariant. Given that $\tau_M$ and $\sigma$ are $G$-invariant, we have
		\begin{align*}
			\tau_M(\pi(a) (s^{-1} \cdot F(sts^{-1}))) &= \tau_M(\pi(s \cdot a) F(sts^{-1})) \\
			&= \sigma((s \cdot a) sts^{-1}) \\
			&= \sigma(at) \\
			&= \tau_M(\pi(a) F(t)),
		\end{align*}
		and so we may apply \Cref{x_t_trace_uniqueness} to conclude that $s^{-1} \cdot F(sts^{-1}) = F(t)$, i.e.\ $s \cdot F(t) = F(sts^{-1})$. This is enough to guarantee $G$-equivariance on the entire domain, as
		$$ F(s \cdot (at)) = F((s \cdot a) sts^{-1}) = \pi(s \cdot a) F(sts^{-1}) = s \cdot (\pi(a) F(t)) = s \cdot F(at). $$
	\end{proof}

	\begin{proof}[Proof of \Cref{crossed_product_traces_reduced}]
		The proof that the given map from set (1) to set (2) is well-defined and injective is analogous to what was done in the proof of \Cref{crossed_product_traces_universal}.
		
		To go from (2) to (3), we first note that $A \rtimes_\lambda R_a(G) = A \rtimes_u R_a(G)$ by amenability of $R_a(G)$, and so there is a trace $\sigma' \in T(A \rtimes_\lambda R_a(G))$ satisfying $\sigma'(a \lambda_t) = \tau_M(\pi(a) x_t)$ by \Cref{crossed_product_traces_universal}. Composing with the canonical conditional expectation $E_{R_a(G)} : A \rtimes_\lambda G \to A \rtimes_\lambda R_a(G)$, which maps $a \lambda_t$ to itself if $t \in R_a(G)$ and zero otherwise, gives us a state $\sigma \defeq \sigma' \circ E_{R_a(G)} \in S(A \rtimes_\lambda G)$. It remains to check that this is indeed still a trace on $A \rtimes_\lambda G$. Note that for $s,t \in G$, we have $st \in R_a(G)$ if and only if $ts \in R_a(G)$ by normality of $R_a(G)$. Hence, if $st \notin R_a(G)$, then
		$$ \sigma(a\lambda_s b\lambda_t) = \sigma(a (s \cdot b) \lambda_{st}) = 0 = \sigma(b (t \cdot a) \lambda_{ts}) = \sigma(b \lambda_t a \lambda_s). $$
		The case of $st \in R_a(G)$ is identical to what was done in the proof of \Cref{x_t_trace_existence}.
		
		Finally, we again wish to show that the composition of the maps from (1) to (2) and (2) to (3) is surjective, i.e.\ given $\sigma \in T(A \rtimes_\lambda G)$ with $\sigma|_A = \tau$, there exists some weak expectation $F : A \rtimes_u R_a(G) \to M$ for $(A,\tau,G,R_a(G))$ satisfying $\sigma = \tau_M \circ F \circ E_{R_a(G)}$. (This last composition makes sense, as $A \rtimes_u R_a(G) = A \rtimes_\lambda R_a(G)$ by amenability). Letting $\rho : A \rtimes_u G \to A \rtimes_\lambda G$ be the canonical *-homomorphism, we note that $\sigma \circ \rho \in T(A \rtimes_u G)$, and so there is some weak expectation $F' : A \rtimes_u G \to M$ for $(A,\tau,G)$ satisfying $\sigma \circ \rho = \tau_M \circ F'$ by \Cref{crossed_product_traces_universal}. Observe that we canonically have $A \rtimes_u R_a(G) \subseteq A \rtimes_u G$---this is because the following composition of canonical maps yields the identity map:
		$$ A \rtimes_u R_a(G) \to A \rtimes_u G \to A \rtimes_\lambda G \to A \rtimes_\lambda R_a(G) = A \rtimes_u R_a(G) $$		
		We claim that $F \defeq F'|_{A \rtimes_u R_a(G)}$ is the map we are looking for. This follows from \cite[Theorem~5.2]{bryder_kennedy_twisted_crossed_products}, which says that $\sigma = \sigma \circ E_{R_a(G)}$.
	\end{proof}

	\begin{proof}[Proof of \Cref{properly_outer_implies_unique_trace_property}]
		This follows immediately from \Cref{crossed_product_traces_universal,crossed_product_traces_reduced}, and the fact that properly outer and freely acting are equivalent (see the review in  \Cref{subsection_properly_outer_automorphisms}).
	\end{proof}
	
	\begin{proof}[Proof of \Cref{unique_trace_property_almost_inner_actions}]
		This follows immediately from \Cref{crossed_product_traces_universal,crossed_product_traces_reduced}, together with the correspondence given in \Cref{jointly_properly_outer_jointly_freely_acting_equivalence}.
	\end{proof}

	\begin{proof}[Proof of \Cref{unique_trace_property_FC_groups}]
		If the action is jointly almost properly outer relative to $R_a(G)$, then \Cref{unique_trace_property_almost_inner_actions} tells us that $\tau$ has unique tracial extension. Conversely, assume the action is partially almost inner relative to $R_a(G)$ with respect to some nontrivial $\set{(p_t, u_t)}_{t \in R_a(G)}$. Pick $t_0 \neq e$ such that $p_{t_0} \neq 0$, and let $C$ denote the conjugacy class of $t_0$ in $G$. Now define
		$$ 
		v_t
		=
		\begin{cases*}
			1 & if $t = e$ \\
			u_t & if $t \in C \cup C^{-1}$ \\
			0 & otherwise
		\end{cases*},
		q_t
		=
		\begin{cases*}
			1 & if $t = e$ \\
			p_t & if $t \in C \cup C^{-1}$ \\
			0 & otherwise
		\end{cases*},
		y_t
		=
		\begin{cases*}
			1 & if $t = e$ \\
			\frac{1}{2 \card{C}} p_t & if $t \in C \cup C^{-1}$ \\
			0 & otherwise
		\end{cases*}
		$$
		Observe that $[v_{st^{-1}} y_{st^{-1}}]_{s,t \in G} = 1 + \sum_{t \in C \cup C^{-1}} \frac{1}{2 \card{C}} u_t \tensor \lambda_t$ is in fact a positive element in $M \tensor_{\min} C^*_\lambda(G)$, as $\sum_{t \in C \cup C^{-1}} \frac{1}{2 \card{C}} u_t \tensor \lambda_t$ is a self-adjoint element of norm at most $1$, and so the action is partially almost inner relative to $R_a(G)$ with respect to the nontrivial and positively compatible set $\set{(q_t, v_t)}_{t \in R_a(G)}$. By \Cref{unique_trace_property_almost_inner_actions}, we are done.
	\end{proof}

	\begin{proof}[Proof of \Cref{unique_trace_property_abelian_groups}]
		If $\tau$ does not have unique tracial extension, then \Cref{unique_trace_property_FC_groups} says that the action is partially almost inner relative to $R_a(G)$ with respect to some nontrivial $\set{(p_t, u_t)}_{t \in R_a(G)}$. Choosing $t \in R_a(G) \setminus \set{e}$ with $p_t \neq 0$ gives us what we want.
		
		Conversely, assume that we do have such a $t \in R_a(G) \setminus \set{e}$, $p \neq 0$ in $M$, and $u \in U(Mp)$. If $t \neq t^{-1}$, then letting $p_t = p_{t^{-1}} = p$, $u_t = u$, $u_{t^{-1}} = u^*$, and $p_s,u_s = 0$ for $s \neq e, t, t^{-1}$ gives us a nontrivial $\set{(p_t, u_t)}_{t \in R_a(G)}$, and so $\tau$ cannot have unique tracial extension by \Cref{unique_trace_property_FC_groups}.
		
		The case of $t = t^{-1}$ requires just a bit more work. Letting $w = u^2$, observe that $\Ad w = \id$, and so $w \in (Z(Mp))^G$, which is a commutative von Neumann algebra. Even in the non-separable setting, every such algebra is isomorphic to $L^\infty(Y,\nu)$ for some locally compact $Y$ and positive Radon measure $\nu$ on $Y$---see, for example, \cite[Chapter~III, Theorem~1.18]{takesaki_theory_of_operator_algebras_I}. (This is the space of all measurable, locally essentially bounded functions from $Y$ to $\C$, modulo agreeing locally almost everywhere).
		Thus, we may choose a unitary $v \in (Z(Mp))^G$ with the property that $v^2 = w^*$. Now letting $p_t = p$, $u_t = uv$, and $p_s,u_s = 0$ for $s \neq e, t$, we obtain a nontrivial $\set{(p_t, u_t)}_{t \in R_a(G)}$ as before, and again $\tau$ cannot have unique tracial extension by \Cref{unique_trace_property_FC_groups}.
	\end{proof}

	The following results describe how the spectrum of a commutative von Neumann algebra breaks up with respect to a periodic automorphism, and are needed for the proof of \Cref{unique_trace_property_integers}. Recall that an extremally disconnected topological space is one where the closure of any open set is open, and that the spectrum of a commutative von Neumann algebra is always extremally disconnected---see, for example, \cite[Chapter~III, Theorem~1.18]{takesaki_theory_of_operator_algebras_I}.
	In terms of notation, $d|n$ will denote ``$d$ divides $n$''.
	
	\begin{lemma}
	\label{extremally_disconnected_periodic_homeomorphism_orbit_size_decomposition}
		Assume $X$ is an extremally disconnected compact Hausdorff space, and $\alpha : X \to X$ is a homeomorphism satisfying $\alpha^n = \id$ for some $n \in \N$. Then $X$ breaks up as $X = \sqcup_{d | n} X_d$, where each $X_d$ is clopen, $\alpha$-invariant, and has the property that every $x \in X_d$ has orbit of size $d$.
	\end{lemma}

	\begin{proof}
		Letting $Y_d = \Fix(\alpha^d)$ for $d | n$, we know that $Y_d$ is clopen by Frol\'ik's theorem---see \cite[Theorem~3.1]{frolik_extremally_disconnected_fixed_points}. This is the set of all points whose orbit size divides $d$. From here, we can let $X_d = Y_d \setminus (\cup_{\substack{m | n \\ m < d}} Y_m)$.
	\end{proof}

	\begin{lemma}
	\label{extremally_disconnected_periodic_homeomorphism_orbit_clopen_transversal}
		Assume $X$ is an extremally disconnected compact Hausdorff space, and $\alpha : X \to X$ is a homeomorphism with the property that every orbit is finite and of the same size $n \in \N$. Then there is a clopen transversal of the orbits, i.e.\ there is some clopen $E \subseteq X$ with the property that $X = \sqcup_{k=0}^{n-1} \alpha^k(E)$.
	\end{lemma}

	\begin{proof}
		We claim that there is at least one nonempty open subset $U \subseteq X$ with the property that all of $U, \alpha(U), \dots, \alpha^{n-1}(U)$ are pairwise disjoint. To see this, choose any $x \in X$, and let $U_k$, $k=0,\dots,n-1$, be pairwise disjoint open sets satisfying $p^k(x) \in U_k$. Now letting
		$$ U = U_0 \cap p^{-1}(U_1) \cap \dots \cap p^{-(n-1)}(U_{n-1}), $$
		we have that $U, p(U), \dots, p^{n-1}(U)$ are all pairwise disjoint.
		
		Given an ascending chain $(U_\lambda)$ of such open sets, the union $\cup U_\lambda$ is still such a set, and so by Zorn's lemma, there is a maximal open set $U$ with this property. We claim that it is in fact clopen. This follows from the following fact: if $V,W \subseteq X$ are open and $V \cap W = \emptyset$, then $V \cap \closure{W} = \emptyset$, and as $\closure{W}$ is open, we have $\closure{V} \cap \closure{W} = \emptyset$.
		
		Finally, we claim that our maximal set $U$ is in fact the set we are looking for, i.e.\ $X = \sqcup_{k=0}^{n-1} \alpha^k(U)$. Assume otherwise, and consider the smaller space $X \setminus \sqcup_{k=0}^{n-1} \alpha^k(U)$ (a clopen, $\alpha$-invariant subset of $X$). Obtaining as before a nonempty open subset $V \subseteq X \setminus \sqcup_{k=0}^{n-1} \alpha^k(U)$ with the property that $V, \alpha(V), \dots, \alpha^{n-1}(V)$ are all pairwise disjoint, the set $U \cup V$ again satisfies this property, contradicting maximality of $U$.
	\end{proof}

	\begin{proof}[Proof of \Cref{unique_trace_property_integers}]
		First, if the action of $\Z$ on $M$ is properly outer, then $\tau$ has unique tracial extension by \Cref{properly_outer_implies_unique_trace_property}. Conversely, assume the action of $\Z$ on $M$ is not properly outer, and let $n \geq 1$ be such that $\alpha^n$ is not properly outer on $M$. Let $p$ be the largest $\alpha^n$-invariant central projection such that $\alpha^n|_{Mp}$ is inner, and fix a unitary $u \in U(Mp)$ implementing this action. Observe that for any $x \in M \alpha(p)$, we have
		$$ \alpha(u) x \alpha(u)^* = \alpha(u \alpha^{-1}(x) u^*) = \alpha \alpha^n \alpha^{-1} x = \alpha^n x. $$
		In other words, $\alpha^n$ is inner on $M \alpha(p)$. By assumption, $\alpha(p) \leq p$. But then
		$$ p = \alpha^n(p) \leq \alpha^{n-1}(p) \leq \dots \leq \alpha(p) \leq p. $$
		This shows $p$ is in fact $\alpha$-invariant. In general, even though the choice of unitary $u \in U(Mp)$ satisfying $\alpha^n = \Ad u$ is not unique, we still cannot guarantee that there is some choice that also satisfies $\alpha(u) = u$---see \cite[Proposition~1.6]{connes_periodic} for an example of this phenomenon on the separable hyperfinite $II_1$ factor. However, we will show that it is always possible to choose an $\alpha$-invariant unitary implementing the action of $\alpha^{n^2}$. From here, \Cref{unique_trace_property_abelian_groups} will apply, giving us the fact that $\tau$ cannot have unique tracial extension.
		
		To simplify notation, we can assume without loss of generality that $p = 1$. Observe that our previous computations above show that $\Ad \alpha(u) = \Ad u$, and so $\alpha(u) = uv$ for some $v \in Z(M)$. Moreover, the the fact that $\alpha^n = \Ad u$ tells us that
		$$ u = \alpha^n(u) = \alpha^{n-1}(uv) = \dots = u v \alpha(v) \dots \alpha^{n-1}(v), $$
		or in other words,
		$$ v \alpha(v) \dots \alpha^{n-1}(v) = 1. $$
		
		Now let $Z(M) = C(X)$. We know that $\alpha$ induces a homeomorphism on $X$, which we will denote by $\alpha_X$. Given that $\alpha^n$ is inner, we know that $\alpha_X^n$ is the identity map. By \Cref{extremally_disconnected_periodic_homeomorphism_orbit_size_decomposition}, we have $X = \sqcup_{d | n} X_d$, where $X_d$ is the set of all $x \in X$ with the size of the $\alpha_X$-orbit being exactly $d$, and furthermore each $X_d$ is clopen and $\alpha_X$-invariant.
		
		We will show for every $d | n$ that there is some central unitary $w_d \in M 1_{X_d}$ with the property that $\alpha(u^n w_d) = u^n w_d$. Again to simplify notation, we may assume without loss of generality that $X = X_d$ for a single $d | n$. Applying \Cref{extremally_disconnected_periodic_homeomorphism_orbit_clopen_transversal}, there exists a clopen transversal $E \subseteq X$ of the orbits of $\alpha_X$. Let $q = 1_E$, and observe that $q, \dots, \alpha^{d-1}(q)$ are pairwise orthogonal projections that sum to $1$. Keeping this in mind, we may decompose $v$ as follows: let $v_k = \alpha^{-k}(v) q \in Mq$ for $k = 0,\dots,d-1$, so that
		$$ v = v_0 + \dots + \alpha^{d-1}(v_{d-1}). $$
		Now define $w$ as follows:
		$$ w = q + \alpha(v_1^n) + \alpha^2(v_1^n v_2^n) + \dots + \alpha^{d-1}(v_1^n \dots v_{d-1}^n), $$
		and note that
		$$ \alpha(w)^* = (v_1^n \dots v_{d-1}^n)^* + \alpha(q) + \alpha^2(v_1^n)^* + \dots + \alpha^{d-1}(v_1^n + \dots v_{d-2}^n)^*, $$
		so that
		$$ w \alpha(w)^* = (v_1^n \dots v_{d-1}^n)^* + \alpha(v_1^n) + \alpha^2(v_2^n) + \dots + \alpha^{d-1}(v_{d-1}^n). $$
		We claim that we in fact have $w \alpha(w)^* = v^n$. Our earlier equality $v \dots \alpha^{n-1}(v) = 1$ gives us $(v_0 \dots v_{d-1})^{n/d} = 1$, and so we obtain the equality $(v_1^n \dots v_{d-1}^n)^* = v_0^n$.
		
		In summary, we have obtained a central unitary $w$ with the property that $\alpha(w) = (v^n)^*w$. Keeping in mind that $\alpha^{n^2} = \Ad u^n = \Ad (u^n w)$, and also that
		$$ \alpha(u^n w) = \alpha(u)^n \alpha(w) = u^n v^n (v^n)^* w = u^n w, $$
		we may apply \Cref{unique_trace_property_abelian_groups} to conclude that $\tau$ cannot have unique tracial extension to $A \rtimes_\lambda \Z$.
	\end{proof}

	\begin{proof}[Proof of \Cref{unique_trace_property_commutative}]
		The proofs for the universal crossed product and reduced crossed product are almost identical. Hence, we only prove the reduced case.
		
		First, assume that the action of $R_a(G)$ on $(X,\mu)$ is essentially free. This is equivalent to the action of $R_a(G)$ on $L^\infty(X,\mu)$ being properly outer. By \Cref{properly_outer_implies_unique_trace_property}, $\mu$ has unique tracial extension.
		
		Now assume that the action of $R_a(G)$ on $(X,\mu)$ is not essentially free, and let $p_t = u_t = 1_{\Fix(t)}$. It is straightforward to check that the action of $G$ on $L^\infty(X,\mu)$ is partially almost inner relative to $R_a(G)$ with respect to the (nontrivial by assumption) $\set{(p_t,u_t)}_{t \in R_a(G)}$. Furthermore, it is clear that the additional assumptions of \Cref{vanishing_obstruction_implies_positively_compatible} are satisfied, and so $\set{(p_t, u_t)}_{t \in R_a(G)}$ is also positively compatible. Again by \Cref{unique_trace_property_almost_inner_actions}, $\mu$ cannot have unique tracial extension.
	\end{proof}

	\section{Examples}
	
	\subsection{A finite-dimensional (counter)example}
	\label{subsection_counterexample_finite}
	
	Here, we give an action of $\Z_2 \times \Z_2$ on $M_2$ such that the action of each $t \in \Z_2 \times \Z_2$ is inner, but the crossed product is isomorphic to $M_4$. Observe that, letting $\tau$ be the unique (hence, automatically invariant) trace on $M_2$, and $\pi : M_2 \to B(H_\tau)$ be the GNS representation, we canonically have $\pi(M_2)'' \isoto M_2$. In particular, the action of $\Z_2 \times \Z_2$ on $\pi(M_2)''$ is not properly outer, but the only invariant trace on $M_2$ extends to a unique trace on the crossed product. Of course, one could also view the crossed product as a von Neumann crossed product, and from this perspective it is a finite factor of type $I$.
	
	This contradicts \cite[Theorem~4.3]{thomsen_trace_bijection}, which claims that if $G$ is countable and abelian and $A$ is unital and separable, then four various conditions are equivalent. In particular, condition (1), which states that $T(A \rtimes_\lambda G)$ and $T_G(A)$ are in canonical bijection, is equivalent to condition (4), which states that for any $\tau \in \boundary_e(T_G(A))$, letting $\pi : A \to B(H_\tau)$ be the GNS representation, the action of $G$ on $\pi(A)''$ is properly outer. The example in this section contradicts (1) $\implies$ (4). The converse, along with the equivalence between (1), (2), and (3), still appear to be correct.
	
	Similarly, this example also contradicts the precursor result \cite[Proposition~11]{bedos_simple_unique_trace}, which again gives an equivalence between three conditions. It, in particular, claims that if $G$ is abelian and acts on a finite factor $N$, then condition (a) stating the von Neumann crossed product $N \closure{\rtimes} G$ is a finite factor is equivalent to condition (c) stating the action is properly outer. Again, (a) $\implies$ (c) is false, but the converse, along with (a) $\iff$ (b), still appear to be correct.
		
	We first present a proof of the following example using purely elementary techniques, and afterwards show how our results apply.
	
	\begin{example}
	\label{example_M2_V4}
		Consider $G = \Z_2 \times \Z_2 = \left<u\right> \times \left<v\right>$, acting on $A = M_2$, where the action $\alpha : G \to \Aut(M_2)$ is given by $\alpha_u = \Ad U$ and $\alpha_v = \Ad V$, where
		$$ U =
		\begin{bmatrix}
			1 & 0 \\
			0 & -1
		\end{bmatrix},
		\quad
		V
		=
		\begin{bmatrix}
			0 & 1 \\
			1 & 0
		\end{bmatrix}. $$
		Then this is a well-defined action, and the crossed product isomorphic to $M_4$, and is hence simple and has unique trace.
	\end{example}
	
	\begin{proof}[Proof using elementary techniques]
		It is easy to check that $\alpha_s$ and $\alpha_t$ are commuting automorphisms, both of order $2$, and so we obtain a group homomorphism $\alpha : G \to \Aut(M_2)$. Given that the crossed product is 16-dimensional, it suffices to show that it has trivial center in order to prove it is isomorphic to $M_4$.
		
		To this end, assume that $\sum_s a_s \lambda_s \in Z(A \rtimes_\lambda G)$. Then given any $t \in G$, we have
		$$ \lambda_t \left(\sum_s a_s \lambda_s\right) = \sum_s (t \cdot a_s) \lambda_{ts} = \sum_s (t \cdot a_s) \lambda_{st}, $$
		while
		$$ \left(\sum_s a_s \lambda_s\right) \lambda_t = \sum_s a_s \lambda_{st}. $$
		This shows each $a_s$ is invariant under the action of each $t \in G$. In other words, $a_s$ commutes with each of the matrices $I$, $U$, $V$, and $UV$. But these matrices are easily seen to span $M_2$, and so $a_s \in Z(M_2) = \C$.
		
		Now letting $b \in M_2$ be arbitrary, we have
		$$ b\left(\sum_s a_s \lambda_s\right) = \sum_s (ba_s) \lambda_s, $$
		while
		$$ \left(\sum_s a_s \lambda_s\right) b = \sum_s (a_s (s \cdot b)) \lambda_s. $$
		If $s \in G$ is such that $a_s \neq 0$, then $b = s \cdot b$ for all $b$. Writing $\alpha_s = \Ad W$, this tells us $W \in Z(M_2) = \C$, so $s = e$.
	\end{proof}

	\begin{proof}[Proof using \Cref{unique_trace_property_abelian_groups}]
		As the crossed product is 16-dimensional, it suffices to prove that it has unique trace. Note that $M_2$ has a unique trace, and the double commutant under its GNS representation is again $M_2$. Assume there is a nontrivial element $t \in \Z_2 \times \Z_2 \setminus \set{e}$, a nontrivial central projection $p \in M_2$, and a unitary $w \in (M_2) p$ with the properties that
		\begin{enumerate}
			\item $s \cdot p = p$ and $s \cdot w = w$ for all $s \in \Z_2 \times \Z_2$.
			\item $t$ acts by $\Ad w$ on $(M_2) p$.
		\end{enumerate}
		Clearly, we must have $p = 1$. Also, as $Z(M_2) = \C$, if there is one unitary $w \in M_2$ implementing the implementing the action of $t$ and satisfying $s \cdot w = w$ for all $s \in \Z_2 \times \Z_2$, then all unitaries implementing this inner action necessarily satisfy this invariance property. But $UV = -VU$, and so $v \cdot U = -U$, $u \cdot V = -V$, and $v \cdot (UV) = -UV$. Thus, the above situation cannot occur, and so by \Cref{unique_trace_property_abelian_groups}, we are done.
	\end{proof}

	\subsection{A finite cyclic group (counter)example}
	\label{subsection_counterexample_finite_cyclic_group}
	
	This section aims to give another counterexample to results cited in \Cref{subsection_counterexample_finite}, but in the case of $G$ being a finite cyclic group instead, and also show how our results apply.
	
	Before proceeding further, we first recall the notion of \textit{separably inheritable} in the sense of Blackadar \cite[Definition~II.8.5.1]{blackadar_operator_algebras}. We say that a property $(P)$ is \textit{separably inheritable} if the following two hold:
	\begin{enumerate}
		\item Whenever $A$ is a C*-algebra satisfying $(P)$, and $B \subseteq A$ is a separable C*-subalgebra, then there is an intermediate C*-algebra $C$ with $B \subseteq C \subseteq A$ with $C$ separable and satisfying $(P)$.
		
		\item Whenever $A_1 \injectsinto A_2 \injectsinto \dots$ is an inductive system of separable C*-algebras, each satisfying $(P)$, with injective connecting maps, then the direct limit $\varinjlim A_n$ also satisfies $(P)$.
	\end{enumerate}
	It is remarked that, in the unital category, the property of having a unique trace is separably inheritable \cite[II.8.5.5]{blackadar_operator_algebras}. The following lemma shows that this works in the unital \textit{equivariant} category as well, and will allow us to construct a counterexample in the separable setting.
	
	\begin{lemma}
	\label{factor_dense_subalgebra_unique_trace}
		Let $G$ be a countable discrete group, and consider the category of unital $G$-C*-algebras.
		\begin{enumerate}
			\item If $A$ is a C*-algebra in this category with a unique trace, and $B$ is any separable C*-subalgebra (not necessarily unital or $G$-invariant), then there is a unital $G$-invariant C*-subalgebra $C \subseteq A$ with $C$ separable and having a unique trace, and also containing $B$.
			
			\item Whenever $A_1 \injectsinto A_2 \injectsinto \dots$ is an inductive system of separable C*-algebras in this category, each with a unique trace, with injective connecting maps, then the direct limit $\varinjlim A_n$ also has a unique trace.
		\end{enumerate}
	\end{lemma}

	\begin{proof}
		The fact that the inductive limit property works is clear. Hence, we prove the intermediate C*-algebra property. Let $A$ and $B$ be as above. We claim that there is always a unital $G$-invariant intermediate C*-subalgebra $\wtilde{B}$ satisfying $B \subseteq \wtilde{B} \subseteq A$. To see this, let $(b_n) \subseteq B$ be a norm-dense sequence, and let
		$$ \wtilde{B} \defeq C^*(1,\setbuilder{g b_n}{n \in \N, g \in G}). $$
		From here, we may construct a sequence of subalgebras of $A$ satisfying
		$$ A_1 \subseteq \wtilde{A}_1 \subseteq A_2 \subseteq \wtilde{A}_2 \subseteq \dots, $$
		where $A_1 = B$, $\wtilde{A}_n$ is defined as before in relation to $A_n$, and $A_{n+1}$ has a unique trace. Then it is clear that
		$$ C \defeq \closure{\bigcup_{n \in \N} A_n} = \closure{\bigcup_{n \in \N} \wtilde{A}_n} $$
		will satisfy the properties we want.
	\end{proof}

	\begin{example}
		There is a separable C*-algebra $A$ and an automorphism $\alpha \in \Aut(A)$ of order $4$ such that the following is true: $A$ has a unique trace $\tau$, and if we let $\pi : A \to B(H_\tau)$ be the GNS representation of $(A,\tau)$ and $M = \pi(A)''$, then:
		\begin{enumerate}
			\item $M$ is the separable hyperfinite $II_1$ factor.
			\item The corresponding action of $\Z_4$ on $M$ is not (properly) outer. In fact, the action of $2 \in \Z_4$ on $A$ is inner (and hence also inner on $M$).
			\item The C*-crossed product $A \rtimes_\lambda \Z_4$ has a unique trace.
			\item The von Neumann crossed product $M \closure{\rtimes} \Z_4$ is a $II_1$ factor.
		\end{enumerate} 
	\end{example}

	\begin{proof}
		Let $R$ be the separable hyperfinite $II_1$ factor. It was shown in \cite[Proposition~1.6]{connes_periodic} that for any $p \in \N$ and any $p$-th root of unity $\gamma$, there is an automorphism $s^\gamma_p \in \Aut(R)$ with the properties that $p$ is the smallest positive integer with $(s^\gamma_p)^p$ being inner, $(s^\gamma_p)^p = \Ad U_\gamma$, where writing $R = \closure{\tensor}_{n=1}^\infty M_p$, we have
		$$ U_\gamma =
		\begin{bmatrix}
			\gamma & & & \\
			& \gamma^2 & & \\
			& & \ddots & \\
			& & & \gamma^p
		\end{bmatrix}
		\tensor
		\left(\tensor_{n=2}^\infty I\right),
		$$
		and moreover, $s^\gamma_p(U_\gamma) = \gamma U_\gamma$. Observe that, as $Z(R) = \C$, then $s^\gamma_p(W) = \gamma W$ for any unitary $W \in R$ satisfying $(s^\gamma_p)^p = \Ad W$. For our purposes, we will let $p=2$ and $\gamma = -1$, and fix an outer automorphism $\alpha \in \Aut(R)$ with $\alpha^2 = \Ad u$ and $\alpha(u) = -u$. Observe that $u = U_2$ as defined above guarantees $u^2 = 1$, and so $\alpha^4 = \id$.
		
		Consider the $\Z_4$ action on $R$ induced by $\alpha$, and observe that by weak*-separability of $R$, we have that $u \in R$ is contained in a norm-separable, weak*-dense C*-subalgebra. By \Cref{factor_dense_subalgebra_unique_trace}, there is a $\Z_4$-invariant, norm-separable, weak*-dense unital C*-subalgebra $A \subseteq R$ containing $u$ and having a unique trace (denote the unique trace on both $R$ and $A$ by $\tau$). We claim that this is the C*-algebra we are looking for.
		
		First, we verify that if $\pi : A \to B(H_\tau)$ is the GNS representation, then we get $\pi(A)'' = R$. Denote the GNS Hilbert spaces of $(A,\tau)$ and $(R,\tau)$ by $L^2(A,\tau)$ and $L^2(R,\tau)$, respectively, and note that $L^2(A,\tau) \subseteq L^2(R,\tau)$. As $A$ is SOT-dense in $R \subseteq B(L^2(R,\tau))$, then if $x \in R$ with $(a_\lambda) \subseteq A$ SOT-convergent to $x$, it is easy to see that $(a_\lambda)$ is also $\norm{\cdot}_2$-convergent to $x$. It follows that $L^2(A,\tau) = L^2(R,\tau)$, and so $\pi(A)'' = R$.
		
		By normality, the unique extension of $\alpha|_A \in \Aut(A)$ to $\pi(A)'' = R$ is again $\alpha$. By construction, the corresponding action of $\Z_4$ on $R$ is not (properly) outer, as $\alpha^2$ is inner (and is in fact inner on $A$ by construction).
		
		We wish to apply \Cref{unique_trace_property_abelian_groups} to conclude that the crossed product $A \rtimes_\lambda \Z_4$ has a unique trace. This follows from our previous computations---the only nontrivial $n \in \Z_4$ that admits a nontrivial central projection $p \in R$ with the properties that $\alpha(p) = p$ and $\alpha^n$ is inner on $Rp$ is $n=2$ (together with $p=1$). However, as we saw earlier, it is impossible to choose a unitary $w$ satisfying both $\alpha^2 = \Ad w$ and $\alpha(w) = w$. This gives us that $\tau$ has unique tracial extension to the crossed product $A \rtimes_\lambda \Z_4$.
		
		From here, we can conclude that the von Neumann crossed product $R \closure{\rtimes} \Z_4$ is still a $II_1$ factor. We know it admits at least one faithful normal trace, namely $\tau \circ E$, where $E : R \closure{\rtimes} \Z_4 \to R$ is the canonical expectation. Given any other normal trace, it necessarily agrees with $\tau \circ E$ on $A \rtimes_\lambda \Z_4$, and by normality and weak*-density therefore agrees with $\tau \circ E$ on all of $R \closure{\rtimes} \Z_4$.
	\end{proof}

	\subsection{Various crossed products with reduced group C*-algebras}
	
	Let $C^*_\lambda(G)$ denote the reduced group C*-algebra, $L(G)$ the group von Neumann algebra, and $\tau_\lambda \in T(L(G))$ the canonical trace. We will, furthermore, denote by $\Char(G)$ the set of all group homomorphisms from $G$ to the circle group $\mathbb{T}$. We say that $\Char(G)$ separates the points of $G$ if for any $s \neq t$ in $G$, there is some $\chi \in \Char(G)$ such that $\chi(s) \neq \chi(t)$. Equivalently, for any $t \neq e$, there is some character $\chi \in \Char(G)$ with $\chi(t) \neq 1$. This definition generalizes to any $H \leq G$ and $K \leq \Char(G)$.
	
	The following facts are likely well-known. In particular, this first proposition is proven in greater generality in \cite[Theorem~5.2]{behncke_automorphisms_of_crossed_products}. We provide quick proofs of them for convenience.
	
	\begin{proposition}
	\label{character_automorphism}
		Every $\chi \in \Char(G)$ induces an automorphism $\alpha_\chi$ on $C^*_\lambda(G)$ and $L(G)$ given by mapping $\lambda_t$ to $\chi(t) \lambda(t)$. If $G$ is ICC, then $\alpha_\chi$ is (properly) outer on $L(G)$ for every $\chi \neq 1$.
	\end{proposition}

	\begin{proof}
		Viewing $C^*_\lambda(G) \subseteq L(G) \subseteq B(\ell^2(G))$, we may define a unitary $U_\chi \in U(\ell^2(G))$ mapping $\delta_t$ to $\chi(t) \delta_t$. From here, we see that
		$$ U_\chi \lambda_s U_\chi^* \delta_t = U_\chi \lambda_s (\conj{\chi(t)} \delta_t) = \conj{\chi(t)} U_\chi \delta_{st} = \conj{\chi(t)} \chi(st) \delta_{st} = \chi(s) \lambda_s \delta_t, $$
		i.e.\ $U_\chi \lambda_s U_\chi^* = \chi(s) \lambda_s$. It follows that $\alpha_\chi \defeq \Ad U_\chi$
		induces the automorphism we want on $C^*_\lambda(G)$ and $L(G)$.
		
		Now assume that $G$ is ICC and that $\chi \in \Char(G)$ satisfies $\alpha_\chi = \Ad u$ for some $u \in L(G)$. Write $u \sim \sum_t \alpha_t \lambda_t$, and observe that
		\begin{align*}
			&u \lambda_s u^* = \chi(s) \lambda_s \\
			\iff &u \lambda_s = \chi(s) \lambda_s u \\
			\iff &\sum_t \alpha_t u_{ts} = \sum_t \chi(s) \alpha_t u_{st} \\
			\iff &\sum_t \alpha_{ts^{-1}} u_t = \sum_t \chi(s) \alpha_{s^{-1}t} u_t \\
			\iff &\alpha_{ts^{-1}} = \chi(s) \alpha_{s^{-1}t} \\
			\iff &\alpha_{sts^{-1}} = \chi(s) \alpha_t
		\end{align*}
		It follows from square-summability of $(\alpha_t)_{t \in G}$ that $\alpha_t = 0$ for $t \neq e$. Hence, $\alpha_\chi = \id$, so $\chi = 1$.
	\end{proof}

	\begin{lemma}
	\label{amenable_radical_product}
		For any groups $G$ and $H$, we have $R_a(G \times H) = R_a(G) \times R_a(H)$.
	\end{lemma}

	\begin{proof}
		Let $\pi_G : G \times H \to G$ and $\pi_H : G \times H \to H$ denote the canonical projections. Observe that $\pi_G(R_a(G \times H))$ is an amenable normal subgroup of $G$, and hence $\pi_G(R_a(G \times H)) \subseteq R_a(G)$. Similarly, $\pi_H(R_a(G \times H)) \subseteq R_a(H)$, and so $R_a(G \times H) \subseteq R_a(G) \times R_a(H)$. But $R_a(G) \times R_a(H)$ is an amenable normal subgroup of $G \times H$, and so we get equality.
	\end{proof}

	It is easy to check that for any $t \in G$ and $\chi \in \Char(G)$, we have that $\Ad \lambda_t$ and $\alpha_\chi$ commute. Thus, for any $H \leq G$ and $K \leq \Char(G)$, we have an action of $H \times K$ on $C^*_\lambda(G)$. This action of course cannot be properly outer on $L(G)$ if $H \neq \set{e}$. However, as this next example shows, $\tau_\lambda$ can still have unique tracial extension to the corresponding reduced crossed product.
	
	\begin{example}
		Assume $G$ is ICC, and let $H \leq G$ and $K \leq \Char(G)$. Then $\tau_\lambda$ has unique tracial extension to $C^*_\lambda(G) \rtimes_\lambda (H \times K)$ if and only if $K$ separates the points of $R_a(H)$.
	\end{example}

	\begin{proof}
		We know that the GNS representation of $(C^*_\lambda(G), \tau_\lambda)$ is the canonical representation $\pi : C^*_\lambda(G) \to B(\ell^2(G))$, and so $\pi(C^*_\lambda(G))'' = L(G)$.
		
		First, assume that $K$ separates the points of $R_a(H)$, and assume the action is partially almost inner relative to $R_a(H) \times K$ with respect to $\set{(p_{t,\chi}, u_{t,\chi})}_{(t,\chi) \in R_a(H) \times K}$ (note that $R_a(H \times K) = R_a(H) \times K$ by \Cref{amenable_radical_product}). Observe that by \Cref{character_automorphism}, in order for $p_{t,\chi} \neq 0$, it must be the case that $\chi = 1$, as the action of $t \in H$ is always inner. Assume $p_{t,e} = 1$ for some nontrivial $t \in H \setminus \set{e}$. Then $u_{t,e} = \gamma \lambda_{t,e}$ for some $\gamma \in \mathbb{T}$. By assumption, there is some $\chi \in K$ with the property that $\chi(t) \neq 1$, and so
		$$ \chi \cdot (\gamma \lambda_t) = \gamma \chi(t) \lambda_t \neq \gamma \lambda_t. $$
		This contradicts the definition of being partially almost inner, and therefore the set $\set{(p_{t,\chi}, u_{t,\chi})}_{(t,\chi) \in R_a(H) \times K}$ is trivial. By \Cref{unique_trace_property_almost_inner_actions}, $\tau_\lambda$ must have unique tracial extension.
		
		Now assume $K$ does not separate the points of $R_a(H)$, and let
		$$ N \defeq \setbuilder{h \in R_a(H)}{\chi(h) = 1 \text{ for all } \chi \in K} \neq \set{e}. $$
		Observe that this is still an amenable normal subgroup of $H$. Now, we will define a partial almost inner action as follows. Given $(t,\chi) \in R_a(H) \times K$, let
		$$
		p_{t,\chi} =
		\begin{cases*}
			1 & if $t \in N$ and $\chi = 1$ \\
			0 & otherwise
		\end{cases*}
		\quad \text{and} \quad
		u_{t,\chi} =
		\begin{cases*}
			\lambda_t & if $t \in N$ and $\chi = 1$ \\
			0 & otherwise
		\end{cases*}.
		$$
		It is straightforward to check that the action of $H \times K$ on $C^*_\lambda(G)$ is indeed partially almost inner relative to $R_a(H) \times K$ with respect to $\set{(p_{t,\chi},u_{t,\chi})}_{(t,\chi) \in R_a(H) \times K}$. Moreover, this is positively compatible by \Cref{vanishing_obstruction_implies_positively_compatible}. Thus, by \Cref{unique_trace_property_almost_inner_actions}, $\tau_\lambda$ cannot have unique tracial extension.
	\end{proof}

	It is worth noting that the above example is not vacuous. For example, we could let $G = \F_2$ with canonical generators $a$ and $b$, $H = \left<a\right>$, and $K$ any subgroup that contains a character mapping $a$ to $e^{2 \pi i \theta}$, where $\theta$ is an irrational number.

	\bibliographystyle{amsalpha}
	\bibliography{crossed_product_traces}
\end{document}